\documentclass[a4paper,10pt]{article}
\usepackage{lineno,hyperref}
\modulolinenumbers[5]
\usepackage{amsmath}
\usepackage{amsfonts}
\usepackage{hyperref}
\usepackage{pgf,tikz}
\usepackage{pifont}
\usepackage{natbib}
\usepackage{amsthm}
\newtheorem{theorem}{Theorem}
\newtheorem{lemma}[theorem]{Lemma}
\newtheorem{corollary}[theorem]{Corollary}
\newtheorem{remark}{Remark}

 \usepackage{mathptmx}      % use Times fonts if available on your TeX system
\usepackage{latexsym}
\usepackage{graphicx}
\usepackage{authblk}
\usepackage{appendix}

\usepackage{algorithm,algorithmicx}
\usepackage{algpseudocode}
\usepackage{placeins}
\usepackage[left=2cm,right=2cm,top=2cm,bottom=2cm]{geometry}

\DeclareMathOperator*{\argmin}{\arg\!\min}

\begin{document}

\title{Weak convergence of particle swarm optimization}

\author[1]{Vianney Bruned}
\author[2]{André Mas}
\author[1]{Sylvain Wlodarczyck}

\affil[1]{Schlumberger Petroleum Services, Montpellier, 34000, France.}
\affil[2]{IMAG, Univ Montpellier, CNRS, Montpellier, France.}

\date{}
\maketitle

\begin{abstract}
Particle swarm optimization algorithm is a stochastic meta-heuristic solving
global optimization problems appreciated for its efficiency and simplicity.
It consists in a swarm of interacting particles exploring a domain and
searching a global optimum. The trajectory of the particles has been
well-studied in a deterministic case. More recently authors shed some light
on the stochastic approach to PSO, considering particles as random variables
and studying them with probabilistic and statistical tools. These works
paved the way to the present article. We focus here on weak convergence, the
kind of stochastic convergence that appears in the Central Limit Theorem. We
obtain three main results related to three different frameworks. These three
settings depend either on the regime of the particles (oscillation or fast
convergence) or on the sampling strategy (along the trajectory or in the
swarm). The main application of these results is the construction of
confidence intervals around the local optimum found by the swarm. The
theorems are illustrated by a simulation study.
\end{abstract}

\bigskip


\raggedbottom

\section{Introduction}

\label{intro} The particle swarm optimization algorithm (PSO), based on
social interactions (behaviors of birds) was introduced in \cite%
{eberhart1995new}. Since then PSO has known a great popularity in many
domains and gave birth to many variants of the original algorithm (see \cite%
{zhang2015comprehensive} for a survey of variants and applications). PSO is
a stochastic meta heuristic solving an optimization problem without any
evaluation of the gradient. The algorithm explores the search space in an
intelligent way thanks to a population of particles interacting with each
other and updating at each step their position and their velocity. The
dynamic of the particles relies on two attractors: their personal best
position (historical best position of the particle denoted $p_{n}^{s}$
below), and the neighborhood best position (corresponding to the social
component of the particles, denoted $g_{n}^{s}$ ). In the dynamic equation,
the attractors are linked with a stochastic process in order to explore the
search space. Algorithm \ref{algo:1} refers to the classical version of PSO
with $S$ particles and $N$ iterations. 
\begin{algorithm}[H]
\begin{algorithmic}
\State Initialize the swarm of $S$ particles with random positions $x_{0}^{s}$ and velocities $v_{0}^{s}$ over the  search space.
\For{$n=1$ to $N$}
\State Evaluate the optimization fitness function for each particle. 
\State Update $p^{s}_n$ (personal best position) and $g^{s}_n$ (neighborhood best position).
\State Change velocity ($v_{n}^{s}$) and position ($x_{n}^{s}$) according to the dynamic equation.
\EndFor
\end{algorithmic}
\caption{Classical PSO}
\label{algo:1}
\end{algorithm}

The convergence and stability analysis of PSO are important matters. In the
literature, there are two kinds of convergence:

\begin{itemize}
\item the convergence of the particles towards a local or global optimum.
This convergence is not obtained with the classical version of PSO. \cite%
{van2010convergence} and \cite{schmitt2015particle} proposed a modified
version of PSO to obtain the convergence.

\item the convergence of each particle to a point (e.g. (\cite{poli2009mean}%
)).
\end{itemize}

If we focus on the convergence of each particle to a point, a prerequisite
is the stability of the trajectory of the particles. In a deterministic
case, \cite{clerc2002particle} dealt with the stability of the particles
with some conditions on the parametrization of PSO. Later, \cite%
{kadirkamanathan2006stability} used the Lyapunov stability theorem to study
the stability. About the convergence of PSO, \cite{van2006study} looked at
the trajectories of the particles and proved that each particle converges to
a stable point (deterministic analysis). Under stagnation hypotheses (no
improvement of the personal and neighborhood best positions), \cite%
{poli2009mean} gives the exact formula of the second moment. More recently, 
\cite{bonyadi2016stability} or \cite{cleghorn2018particle} provided results
for the order-1 and order-2 stabilities with respectively stagnant and
non-stagnant distribution assumptions (both weaker than the stagnation
hypotheses). Since our main results rely on martingale convergence theorems
the other recent work \cite{Xu2018} has to be mentioned as well.

Let us introduce some notations. We consider here a cost function $f:\mathbb{%
R}^{d}\rightarrow\mathbb{R}^{+}$ that should be minimized on a compact set $%
\mathrm{\Omega}$. Consequently the particles evolve in $\mathrm{\Omega}%
\subset\mathbb{R}^{d}$.

Let $x_{n}^{s}\in\mathbb{R}^{d}$ $1\leq s\leq S$ denote the position of
particle number $s$ in the swarm at step $n$. Let $\left( r_{j,n}\right)
_{j=1,...,n\geq1}$ be sequences of independent random vectors in $\mathbb{R}%
^{d}$ whose margins are uniformly distributed over $\left[ 0,1\right] $ and
denote by $\omega,$ $c_{1}$ and $c_{2}$ three positive constants which will
be discussed later. Then the PSO algorithm considered in the sequel is
defined by the two following equations (or dynamic equations, \cite%
{poli2009mean}): 
\begin{equation}
\left\{ \begin{aligned} v_{n+1}^{s}&=\omega \cdot v_{n}^{s}+
c_{1}r_{1,n}\odot\left(p_{n}^{s}-x_{n}^{s}\right)+c_{2}r_{2,n}\odot%
\left(g_{n}^{s}-x_{n}^{s}\right),\\ x_{n+1}^{s}&=x_{n}^{s}+v_{n+1}^{s} \\
\end{aligned}\right.  \label{pso-def}
\end{equation}

where $\odot$ stands for the Hadamard product: 
\[
u\odot v=\left( u_{1}v_{1},...,u_{d}v_{d}\right) 
\]
and $p_{n}^{s}$ (resp. $g_{n}^{s}$) is the best personal position (resp. the
best neighborhood position of the particle $s$) : 
\begin{align*}
p_{n}^{s} &=\argmin_{t \in \lbrace x_{0}^{s},\ldots,
x_{n}^{s}\rbrace}f\left( t\right), \\
g_{n}^{s}&=\argmin_{t \in \lbrace p_{n}^{s^{\prime}}: s^{\prime}\in\mathcal{V%
}(s)\rbrace }f\left( t\right)
\end{align*}
with $\mathcal{V}(s)$ the neighborhood of particle $s$. This neighborhood
depends on the swarm's topology: if the topology is called global (all the
particles communicate between each other) then $g_{n}^{s}=g_{n}=\argmin_{t
\in \lbrace p_{n}^{1},\ldots, p_{n}^{S} \rbrace }f\left( t\right)$ (see \cite%
{lane2008particle}).

Our main objective is to provide (asymptotic) confidence sets for the global
or local optimum of the cost function $f$. If $g=\argmin_{t \in \mathrm{%
\Omega}}f\left(t\right)$ for some domain $\mathrm{\Omega}$, a confidence
region at level $1-\alpha$ (with $\alpha \in \left(0,1\right)$) for $g$ is a random set $\Lambda \subset \mathbb{R}^d$ such that :%
\[
\mathbb{P}\left(g \in \Lambda \right) \geq 1-\alpha.
\]
The set $\Lambda$ depends on
the swarm and is consequently random due to the random evolution of the
particles. The probability symbol above, $\mathbb{P}$, depends on the
distribution of the particles in the swarm. Let us illustrate the use of
confidence interval for a real-valued PSO. Typically the kind of results we
expect is: $g \in \left[m,M\right]$ with a probability larger than, say 99
\%. This does not aim at yielding a precise estimate for $g$ but defines a
``control area" for $g$, as well as a measure of the imprecision and
variability of PSO.

Convergence of the swarm will not be an issue here. In fact we assume that
the personal and global best converge : see assumptions $\mathbf{A_2}$ and $%
\mathbf{B_2}$ below. We are interested in the ``step after" : localizing the
limit of the particles with high probability, whatever their initialization
and trajectories.

Formally, confidence set estimation forces us to inspect order two terms
(i.e. the rate of convergence), typically convergence of the empirical
variance. The word \textit{asymptotic} just means that the sample size
increases to infinity.

\bigskip

The outline of the paper is the following. In the next section the three
main results are introduced. They are all related to weak convergence of
random variables and vectors (see \cite{billing} for a classical monograph)
and obtained under three different sets of assumptions.

The two first consider the trajectory of single particles. The sample
consists in the path of a fixed particle. We show that two different regimes
should be investigated depending on the limiting behavior of $p_n$ and $g_n$%
. Briefly speaking : if the limits of $p_n$ and $g_n$ are distinct, the
particles oscillate between them (which is a well-known characteristics of
PSO), if the limits of $p_n$ and $g_n$ coincide, then particles converge at
a fast, exponential, rate.

In the oscillating case a classical Central Limit Theorem is obtained
relying essentially on martingale difference techniques. In the
non-oscillating situation, the particle converges quickly and we have to use
random matrices products to obtain a non-standard CLT. As by-products of
these two subsections we will retrieve confidence sets of the form $%
\Lambda\left(x_1^s \ldots x_n^s\right)$, depending on the $n$ positions of
each particle $x^s$.

The third result states another classical CLT. The sample consists here in
the whole swarm. This time the confidence set is of the form $%
\Lambda\left(x_n^1 \ldots x_n^S\right)$, depending on the $S$ particles of
the swarm when the iteration step is fixed at $n$. A numerical study and
simulations are performed in a Python environment. A discussion follows. The
derivations of the main theorems are collected in the last section and in
the Appendix.

\section{Main results}

The usual euclidean norm and associated inner product for vectors in $%
\mathbb{R}^{d}$ are denoted respectively $\left\Vert \cdot\right\Vert $ and $%
\left\langle \cdot,\cdot\right\rangle $. If $X$ is a random vector with null
expectation then $\mathbb{E}\left( X\otimes X\right)=\mathbb{E}%
\left(XX^{t}\right)$ is the covariance matrix of $X$. The covariance matrix
is crucial since it determines the limiting Gaussian distribution in the
Central Limit Theorem. We will need two kinds of stochastic convergence in
the sequel. convergence in probability of $X_{n}$ to $X$ is denoted $%
X_{n}\rightarrow_{ \mathbb{P}}X$. The arrow $\hookrightarrow$ stands for
convergence in distribution (weak convergence).

Except in section \ref{sfs} we consider a single particle in order to
alleviate notations. We drop the particle index so that $x_{n}^{s}=x_{n}$, $%
p_{n}^{s}=p_{n}$ and $g_{n}^{s}=g_{n}$.

At last we take for granted that particles are warm, reached an area of the
domain were they fluctuate without exiting (condition $\mathbf{A}_{1}$
below).

\subsection{First case: oscillatory ($p\neq g$)}

%\paragraph*{}

Denote $\xi_{n}=\max\left\{  \left\vert p_{n}-p\right\vert ,\left\vert
g_{n}-g\right\vert \right\}  $ and :%
\[
\theta=\frac{c_{1}p+c_{2}g}{c_{1}+c_{2}},\quad c=\frac{c_{1}+c_{2}}{2}.
\]

The following assumptions are required and discussed after the statement of
Theorem \ref{TH1}.

\noindent $\mathbf{A}_{1}:$ For all $n$, $x_{n} \in \mathrm{\Omega}$ where $%
\mathrm{\Omega}$ is a compact subset of $\mathbb{R}^{d}$.

\noindent $\mathbf{A}_{2} :\sqrt{N}\xi_{N} \overset{\mathbb{P}}{\to} 0$ when $N \to +\infty$.

\noindent $\mathbf{A}_{3}:$ The inequality below connects $c_1$, $c_1$ and $\omega$ : 
\[
2c\frac{1-\omega}{1+\omega}\left(  1+\omega-\frac{c}{2}\right)>\frac{c_{1}^{2}+c_{2}^{2}}{12}
\]

Before stating the Theorem we need a last notation. Let $\delta=\left(
\delta_{1},...,\delta_{d}\right) \in\mathbb{R}^{d}$. The notation $\mathrm{%
diag}\left( \delta\right) $ stands for the diagonal $d\times d$ matrix with
entries $\delta_{1},...,\delta_{d}$ and $\delta^{\odot2}$ is the vector in $%
\mathbb{R}^{d}$ defined by $\delta^{\odot2}=\left(
\delta_{1}^{2},...,\delta_{d}^{2}\right)$.

\begin{theorem}
\label{TH1}Set $\mathfrak{L}=c\left(  \frac{1-\omega}{1+\omega}\right)  \left(
1+\omega-\frac{c}{2}\right) $ and 
\[
\mathfrak{C}=\frac{1}{24}\frac{c_1 c_2}{c^2}\frac{\mathfrak{L}}{\mathfrak{L}-\frac{c_{1}^2+c_{2}^2}{24}}.
\]

Denote finally $\Gamma=\mathfrak{C}\cdot\mathrm{diag}\left(  p-g\right)^{\odot2}$ then :
\[
\sqrt{N}\left(  \frac{1}{N}\sum_{n=1}^{N}x_{n}-\theta\right)  \hookrightarrow
\mathcal{N}\left(  0,\Gamma\right)
\]
where $\mathcal{N}\left(  0,\Gamma\right)  $ denotes the Gaussian centered
random vector of $\mathbb{R}^{d}$ with covariance matrix $\Gamma$.
\end{theorem}

\noindent \textbf{Discussion of the Assumptions:}

\noindent We avoid here the assumption of stagnation: the personal and local
best are not supposed to be constant but they oscillate around their
expectation. The convergence occurs at a rate ensuring that neither $g_{n}$
nor $p_{n}$ are involved in the weak convergence of the particles $x_{n}$.
Condition $\mathbf{A}_{2}$ is specific of what we intend by a convergent
PSO. It ensures that $p_n$ or $g_n$ have no impact on the weak convergence
behavior of the particles. With other words Assumption $\mathbf{A}_{2}$
requires that the oscillations of $p_{n}$ and $g_{n}$ around their
expectations are negligible. We tried here to model the stagnation
phenomenon which consists in sequence of iterations during which $g_{n}$
(resp. $p_{n}$) remain constant for $n$ in $%
\left[ \underline{N}, \overline{N}\right]$ where $\underline{N}<\overline{N}$ are two integers.

Note that assumption $\mathbf{A}_{3}$ is exactly the condition found
in \cite{poli2009mean} (see the last paragraph of section III) for defining
order 2 stability. This condition may be extended to the case when $c_1 \neq
c_2$, see \cite{cleghorn2018particle} and references therein. At last $%
\mathbf{A}_{3}$ holds for the classical calibration appearing in \cite%
{clerc2002particle} (constriction constraints) with $c=1.496172$ and $%
\omega=0.72984$. 
\begin{figure}[]
\centering
\includegraphics[width=0.8\textwidth]{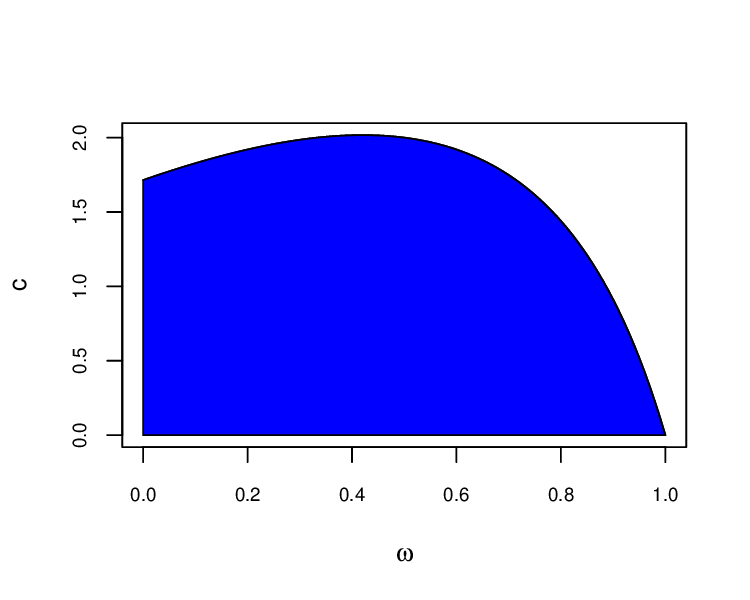}  
\caption{Display of constraint $\mathbf{A}_{3}$ in the plane $\left(\protect%
\omega,c\right)$.}
\label{fig:assumptionA3}
\end{figure}

The next Corollary provides finally two kinds of by-products : asymptotic
confidence intervals, in $\mathbb{R}$, for the coordinates of $\theta$ and
confidence regions for $\theta$ in $\mathbb{R}^d$. Let $\alpha \in \left[0,1%
\right]$, $q_{\alpha}$ be the $1-\alpha$ quantile of the standard Gaussian
distribution and $\chi_{1-\alpha}^{\left(d\right)}$ be the $1-\alpha$
quantile of the Chi-square distribution with $d$ degrees of freedom. Set
also $\theta=\left(\theta_1,\ldots,\theta_d\right)$ and $\bar{x}_N=\left(%
\bar{x}_{N,1},\ldots,\bar{x}_{N,d}\right)$.

\begin{corollary}
Pick $\ell \in \left\{1,\cdots,d\right\}$. An asymptotic (in $N$) confidence interval at level $1-\alpha$ for $\theta_\ell$ is directly derived from Theorem \ref{TH1} :
\[
\mathcal{I}_{1-\alpha}\left(\theta_\ell\right)=\left[\bar{x}_{N,\ell}-s_\ell\left(N,\alpha\right);\bar{x}_{N,\ell}+s_\ell\left(N,\alpha\right)\right]
\]
with $s_\ell\left(N,\alpha\right)=\left|p_\ell-g_\ell\right|\sqrt{\frac{\mathfrak{C}}{N}}q_{1-\frac{\alpha}{2}}$ . \\
An asymptotic confidence (in $N$) region at level $1-\alpha$ for the vector $\theta=\left( p+g\right)/2$ is :
\begin{align*}
\Lambda_{d}\left(1-\alpha\right)=
\left\{t \in \mathbb{R}^d : N\left\|\Gamma^{-1/2}\left(t-\bar{x}_N\right)\right\|^2
\leq \chi_{1-\alpha}^{\left(d\right)} \right\} \\
=\left\{t=\left(t_1,...,t_n\right): \sum_{\ell=1}^{d}\left(\frac{t_\ell-\bar{x}_{N,\ell}}{p_\ell-g_\ell}\right)^2
\leq \frac{\mathfrak{C}}{N}\chi_{1-\alpha}^{\left(d\right)} \right\}.
\end{align*}

\end{corollary}

We note however that the vector $\theta$ may not be of crucial interest for
the initial optimization problem conversely to $g$. This point will be
discussed in the last section.

%First, notice that $\frac{1}{\sqrt{N}}\sum_{n=1}^{N}t_{n}\rightarrow_{ \mathbb{P}}0$ does not imply $\lim_{n\rightarrow+\infty}\mathbb{E}\left\Vert t_{n}\right\Vert =0$ nor the converse. Take for instance $t_{n}=\left(
%-1\right)  ^{n}t_{0}$ with $\mathbb{E}t_{0}\neq0$ then $\frac{1}{ \sqrt{N}
%}\sum_{n=1}^{N}t_{n}\rightarrow_{\mathbb{P}}0$ whereas $\mathbb{E} \left\Vert
%t_{n}\right\Vert =\mathbb{E}\left\Vert t_{0}\right\Vert $. Conversely take
%$t_{n}=t_{0}/\log\left(  n\right)  $ $\lim_{n\rightarrow+\infty}
%\mathbb{E}\left\Vert t_{n}\right\Vert =0$ but $\frac{1}{\sqrt{N}}\sum
%_{n=1}^{N}t_{n}=\frac{t_{0}}{\sqrt{N}}\sum_{n=1}^{N}\log^{-1}\left(  n\right)
%$ cannot converge in probability to 0. 

\subsection{Second case: non-oscillatory and stagnant ($p=g$)}

In this section we study again a single particle and suppose once and for
all that $x_{n}\in\mathbb{R}$. We assume throughout this subsection that the
particle is under stagnation that is $p_n=p$ for $n$ sufficiently large (see
assumption $\mathbf{B_2}$ below). This assumption is strong but a more
general framework leads to theoretical developments out of our scope.
Starting from Equation (\ref{pso-def}), the PSO equation becomes this time: 
\[
x_{n+1}=\left( 1+\omega\right) x_{n}-\omega x_{n-1}+c\left( r_{1,n}
+r_{2,n}\right) \left( p-x_{n}\right). 
\]

%The appendix provides additional material to understand why we confine ourselves to stagnation.
Change the centering and consider $x_{n}-p=y_{n}$. The previous equation
becomes~: 
\begin{equation}  \label{chain:y}
y_{n+1}=\left( 1+\omega-c+c\varepsilon_{n}\right) y_{n}-\omega y_{n-1},
\end{equation}

where $\varepsilon_{n}$ is the sum of two independent random variables with $%
\mathcal{U}\left[ -1/2;1/2\right] $ distribution.

Assuming that for all $n$ $y_{n}\neq0$, we have then : 
\begin{equation}  \label{chain:X}
\frac{y_{n+1}}{y_{n}} =\left( 1+\omega-c+c\varepsilon_{n}\right) -\omega%
\frac{y_{n-1}}{y_{n}}
\end{equation}

It is plain that $y_{n+1}/y_{n}$ defines a Markov chain (more precisely : a
non-linear auto-regressive process) which will play a crucial role in the
forthcoming results. It is shown in the proof section that $y_{n+1}/y_{n}$
is Harris recurrent and has consequently a stationary distribution denoted $%
\pi$. The definition of Harris recurrence needed is given for instance in \cite%
{meyn2012markov}, beginning of Chapter 9. Take $Z_n$ a copy of $y_{n}/y_{n-1}
$ with $Z_0$ a realization of $\pi$. Then define : 
\begin{align*}
\mu_{x} & =\mathbb{E}_{\pi}\log\left\vert Z_0\right\vert , \\
\sigma_{x}^{2} & =\mathrm{Var}_{\pi}\left( \log\left\vert Z_{0}\right\vert
\right) +2\sum_{k=1}^{+\infty}\mathrm{Cov}_{\pi}\left( \log\left\vert
Z_{0}\right\vert ,\log\left\vert Z_{k}\right\vert \right).
\end{align*}

We are ready to introduce a new set of assumptions. \smallskip

\noindent $\mathbf{B_1}: 1+\omega -c< \omega/c < (1+c)/4$.

\noindent $\mathbf{B_2}:$ For sufficiently large $n$ $g_{n}=p_{n}=p=g$ is
constant.

Before stating next Theorem notice that $\sum_{n=1}^{N}\log\left\vert \frac{%
y_{n}}{y_{n-1}}\right\vert =\log\left\vert y_{N}\right\vert -\log\left\vert
y_{0}\right\vert $.

\smallskip

\begin{figure}[]
\centering
\includegraphics[width=0.8\textwidth]{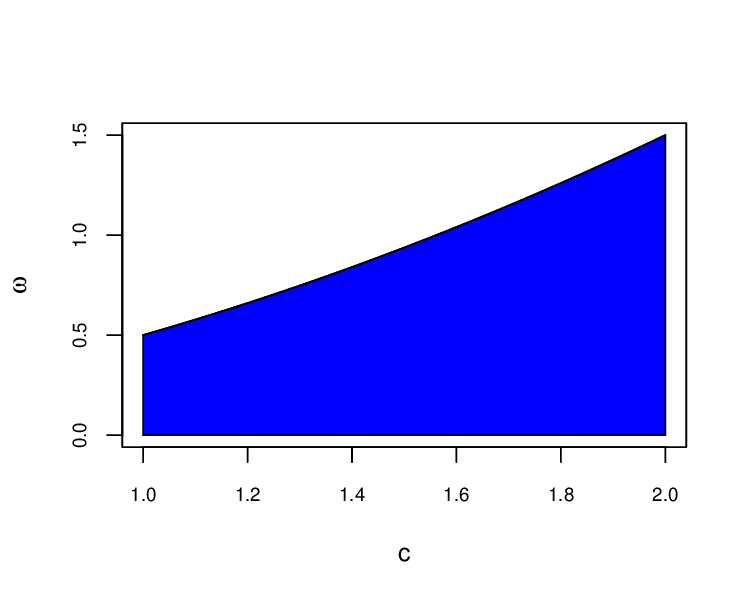}  
\caption{Display of constraint $\mathbf{B}_{1}$ in the plane $\left(c,%
\protect\omega\right)$.}
\label{fig:assumptionB2}
\end{figure}

\begin{theorem}
Let $\omega\in(0,1)$, $c>1$, when $\mathbf{B_{1-2}}$ hold, then :
\[
\frac{1}{\sqrt{N}}\left(  \log\left\vert x_{N}-g_{N}\right\vert -N\mu_{x}\right)
\hookrightarrow\mathcal{N}_1\left(  0,\sigma_{x}^{2}\right)
\]
when $N$ tends to infinity.
\label{theo:equal}
\end{theorem}

\begin{remark}
The theorem above is not a Central Limit Theorem for $x_N$. It is derived thanks to a CLT but it shows that the asymptotic distribution of $\left\vert x_{N}-g_{N}\right\vert$ is asymptotically log-normal with approximate parameters $N\mu_{x}$ and $N\sigma_{x}^{2}$.
\end{remark}

\begin{remark}
The mean and variance $\mu_{x}$ and $\sigma_{x}^{2}$ are usually unknown but
may be approximated numerically. We refer to the simulation section for more details.
\end{remark}

\begin{corollary}
If $p_{n}=p$ for all $n$ (pure stagnation) an asymptotic non convex confidence region for $g$ at level $1-\alpha$ denoted $\Lambda_{1-\alpha}$ below may be derived from the
preceding Theorem:
\begin{align*}
\Lambda_{1-\alpha}\left(g\right)  & =\Lambda_{1-\alpha}^{+}\cup\Lambda_{1-\alpha}^{-}\\
\Lambda_{1-\alpha}^{+}  & =\left[  x_{n}+\exp\left(  \mu_{x}+\frac
{\sigma_{x}}{\sqrt{n}}q_{\alpha/2}\right)  ,x_{n}+\exp\left(  \mu
_{x}-\frac{\sigma_{x}}{\sqrt{n}}q_{1-\alpha/2}\right)  \right]  \\
\Lambda_{1-\alpha}^{-}  & =\left[  x_{n}-\exp\left(  \mu_{x}+\frac
{\sigma_{x}}{\sqrt{n}}q_{1-\alpha/2}\right)  ,x_{n}-\exp\left(  \mu
_{x}-\frac{\sigma_{x}}{\sqrt{n}}q_{\alpha/2}\right)  \right]
\end{align*}

\end{corollary}

\begin{remark}
Under matrix form the equation (\ref{chain:y}) is purely linear but driven by a random matrix:
\begin{align}
\left(
\begin{array}
[c]{c}
y_{n+1}\\
y_{n}
\end{array}
\right)   &  =\mathbf{y}_{n+1} =\left[
\begin{array}
[c]{cc}
1+\omega-c\left(  1+\varepsilon_{n}\right)  & -\omega\\
1 & 0
\end{array}
\right]  \left(
\begin{array}
[c]{c}
y_{n}\\
y_{n-1}
\end{array}
\right) \label{non-osc-pso}\\
\mathbf{y}_{n+1}  &  \mathbf{=}\mathbf{S}_{n+1}\mathbf{y}_{n},\quad 
\mathbf{y}_{n}=\mathbf{S}_{n}\mathbf{S}_{n-1}...\mathbf{S}_{2}\mathbf{y}
_{1}=\mathbf{T}_{n}\mathbf{y}_{1},\nonumber
\end{align}
with $\mathbf{T}_{n}=\mathbf{S}_{n}\mathbf{S}_{n-1}...\mathbf{S}_{2}$. It is
plain here that a classical Central Limit Theorem cannot hold for the sequence
$\left(  y_{n}\right)  _{n\in\mathbb{N}}$.  In the proofs we turn to asymptotic theory for the product of random matrices. We refer to the historical references: \cite{furstenberg1960} and \cite{berger1984central} who proved Central Limit Theorems for the regularity index of the product of i.i.d random matrices. Later \cite{hennion1997limit} generalized their results. But the assumptions of (almost surely) positive entries is common to all these papers. Other authors obtain similar results under different sets of assumptions (see \cite{le1982theoremes}, \cite{benoist2016}, and references therein), typically revolving around characterization of the semi-group spanned by the distribution of $\mathbf{S_n}$. These assumptions are uneasy to check here and we carried out to a direct approach with Markov chain fundamental tools.
\end{remark}

\subsection{The swarm at a fixed step}

\label{sfs}

In this section we change our viewpoint. Instead of considering a single
particle and sampling along its trajectory we will take advantage of the
whole swarm but at a fixed and common iteration step. Our aim here is to
localize the minimum of the cost function based on $\left(x_n^1,\ldots,x_n^S%
\right)$. This time the particle index $s$ varies up to $S$ the swarm size,
whereas the index $n$ is fixed. In this subsection we assume that $S\uparrow
+\infty$ and asymptotic is with respect to $S$. We do not drop $n$ in order
to see how the iteration steps influence the results. We still address only
the case $x_{n}^s\in \mathbb{R}$ even if our results may be
straigthforwardly generalized to $x_{n}^s\in \mathbb{R}^{d}$. We provide
below a Central Limit Theorem suited to the case when the number of
particles in the swarm becomes large. In order to clarify the method, we
assume that for all particles $x_{n}^{i}$ in the swarm $p_{n}^{i}=g_{n}=p$.
In other words, no local minimum stands in the domain $\mathrm{\Omega}$,
which implies additional smoothness or convexity assumptions on the cost
function $f$. This may be possible by a preliminary screening of the search
space. Indeed a first (or several) run(s) of preliminary PSO(s) on the whole
domain identifies an area where a single optimum lies. Then a new PSO is
launched with initial values close to this optimum and with parameters
ensuring that most of the particles will stay in the identified area.

So we are given $\left(x_{n}^{1},...,x_{n}^{S}\right)$ where $S$ is the
sample size. Basically, the framework is the same as in the non oscillatory
case studied above for a single particle. From (\ref{non-osc-pso}) we get
with $\mathbf{y}_{n}^{i}=x_{n}^{i}-p$ : 
\begin{eqnarray*}
\mathbf{y}_{n}^{i} &=&\mathbf{T}_{n}^{i}\mathbf{y}_{1}^{i}, \\
\mathbf{T}_{n}^{i} &=&\Pi _{j=2}^{n}\mathbf{S}_{j},\quad \mathbf{S}_{j}=%
\left[ 
\begin{array}{cc}
1+\omega -c\left( 1+\varepsilon _{j}^{i}\right) & -\omega \\ 
1 & 0%
\end{array}
\right].
\end{eqnarray*}
Assume that the domain $\mathrm{\Omega} $ contains $0$ and that for all $s$ $%
\left( x_{0}^{i},x_{1}^{i}\right) _{i\leq S}$ are independent, identically
distributed and centered then from the decomposition above, for all $n$ and $%
s$,  $\mathbb{E}\mathbf{y}_{n}^{i}=0$ and the $\left( \mathbf{y}%
_{n}^{i}\right) _{1\leq i\leq S}$ are i.i.d too.

\smallskip

The assumptions we need to derive Theorem \ref{prop3} below are :

$\mathbf{C}_{1}:$ The operational domain $\mathrm{\Omega} $ contains $0$
(and is ideally a symmetric set).

$\mathbf{C}_{2}:$ The couples $\left( x_{0}^{i},x_{1}^{i}\right) _{i\leq S}$
are i.i.d. and centered.

$\mathbf{C}_{3}:$ For all $i$ in $\left\{ 1,...,S\right\} $ $%
p_{n}^{i}=g_{n}=p$.

\smallskip

When $S$ is large the following Theorem may be of interest and is a simple
consequence of the i.i.d. CLT.

\begin{theorem}
Under assumptions $\mathbf{C}_{1-3}$ a Central Limit Theorem holds when $S$ the number of
particles in the swarm tends to $+\infty$ :
\begin{equation*}
\frac{1}{\sqrt{S}}\sum_{i=1}^{S}\left( x_{n}^{i}-g_{n}\right) \underset{%
S\rightarrow +\infty }{\hookrightarrow }\mathcal{N}\left( 0,\sigma_{n}^{2}\right),
\end{equation*}
where $\sigma _{n}^{2}=\mathbb{E}\left( x_{n}^{1}-g_{n}\right) ^{2}$ is estimated consistently by :
\begin{equation*}
\widehat{\sigma }_{n}^{2}=\frac{1}{S}\sum_{i=1}^{S}\left(x_{n}^{i}-g_{n}\right) ^{2}.
\end{equation*}
\label{prop3}
\end{theorem}

\begin{remark}
The convergence of $\widehat{\sigma }_{n}^{2}$ to $\sigma _{n}^{2}$ is a
straightforward consequence of the weak and strong laws of large numbers.
\end{remark}

\noindent Denote $\overline{x}^{S}=\left( 1/S\right) \sum_{i=1}^{S}x_{n}^{i}$%
. The Theorem above paves the way towards an asymptotic confidence interval.

\begin{corollary}
An asymptotic confidence interval at level $1-\alpha$ for $g$ is :
\begin{equation*}
\Lambda_n\left(g\right)=\left[ \overline{x}^{S}-\frac{\widehat{\sigma }_{n}}{\sqrt{S}}q_{1-\alpha
/2},\overline{x}^{S}+\frac{\widehat{\sigma }_{n}}{\sqrt{S}}q_{1-\alpha /2}%
\right].
\end{equation*}
\end{corollary}

\section{Simulation and numerical results}

All along this section we take $c_1=c_2=c$.

\begin{figure}[]
\centering
\includegraphics[width=0.6\textwidth]{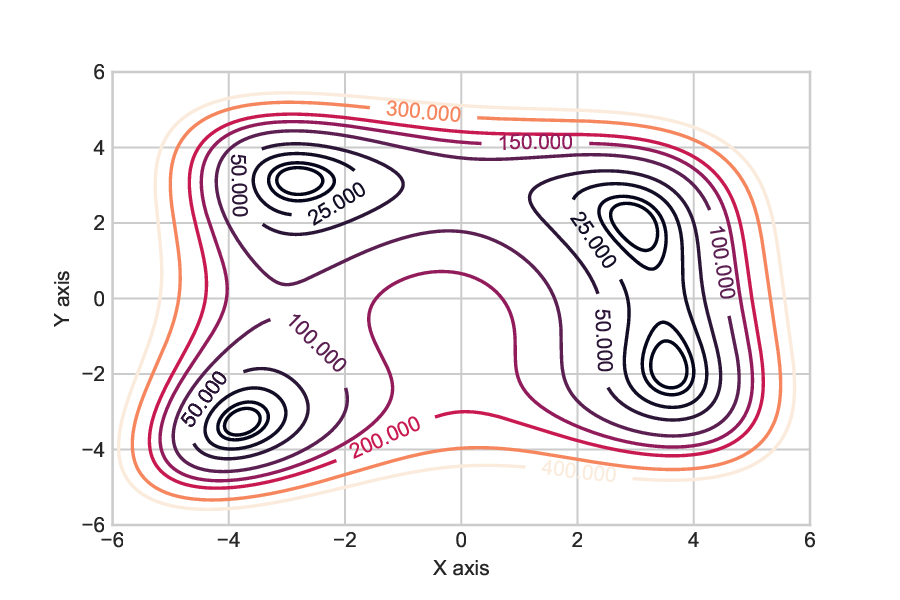}  
\caption{Contour of the Himmelblau's function in the space $\left[-6,6\right]%
^{2}$. }
\label{fig:himmelblau_contour}
\end{figure}

The Himmelblau's function is chosen as example for our experiments. It is a
2 dimensional function with four local optima in $\left[-10,10\right]^{2}$
defined by : $f(x,y)=\left(x^2+y-11\right)^2+\left(x+y^2-7\right)^2$. Figure %
\ref{fig:himmelblau_contour} illustrates the contour of this function.

With the Himmelblau's function, we can observe the two different behaviors
of the particles: oscillatory and non-oscillatory. The Himmelblau's function
is positive and has four global minima in: $(3,2)$, $(-2.81,3.13)$, $%
(-3.77,-3.28)$, $(3.58,-1.84)$ where $f(x,y)=0$. We use a ring topology (for
a quick review of the different topologies of PSO see \cite{lane2008particle}%
) for the algorithm in order to have both oscillating and non-oscillating
particles. The latter converge quickly. The former go on running between two
groups of particles converging to two distinct local optima.

\subsection{Oscillatory case}

We select particles oscillating between $(3.58,-1.84)$ and $(3,2)$, these
values could be both their personal best position or their neighborhood best
position. In this case, the convergence of the $p_n$ and $g_n$ to $%
(3.58,-1.84)$ or $(3,2)$ satisfies the conditions of Theorem \ref{TH1}. We have to verify that the Gaussian asymptotic behavior of $%
H^{s}_{1}(N)=\sqrt{N}\left( \frac{1}{N}\sum_{n=1}^{N}x^{s}_{n}-\frac{p+g}{2}%
\right)$ for each $s$ oscillating particle.

We launch PSO with a population of 200 particles and with 2000 iterations, $%
\omega=0.72984$ and $c=1.496172$. A ring topology was used to ensure the
presence of oscillating particles. A particle is said oscillating if between
the 500th and the 2000th iteration, Assumptions $\mathbf{A}_{1-3}$ holds.

\begin{figure}[]
\centering
\begin{minipage}[t]{.45\linewidth}
		\includegraphics[width=\textwidth]{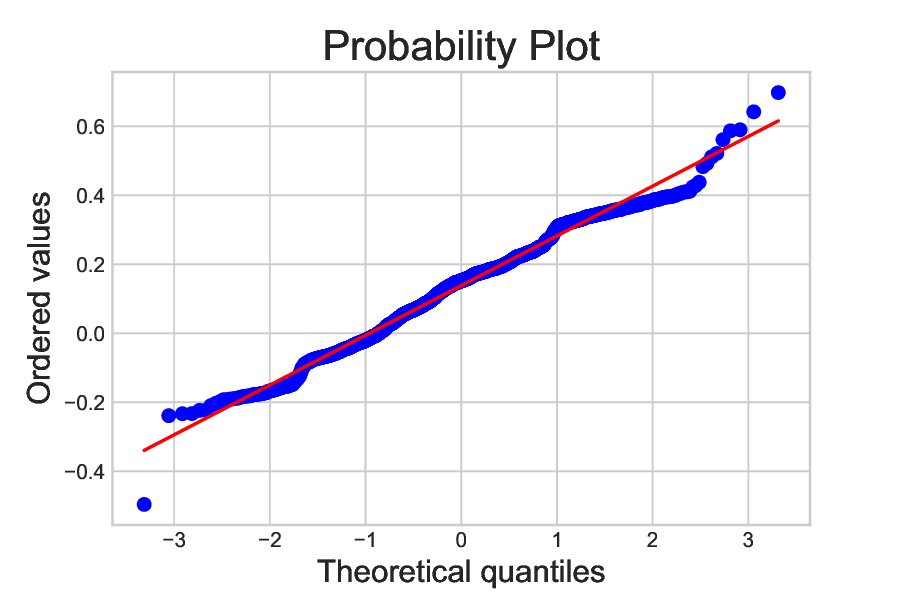}
	\caption{Normal probability plot of $H^{s}_{1}(N)$ on the first coordinate.}
	\label{fig:droite_henry_x}
   \end{minipage}\hfill  
\begin{minipage}[t]{.45\linewidth}
		\includegraphics[width=\textwidth]{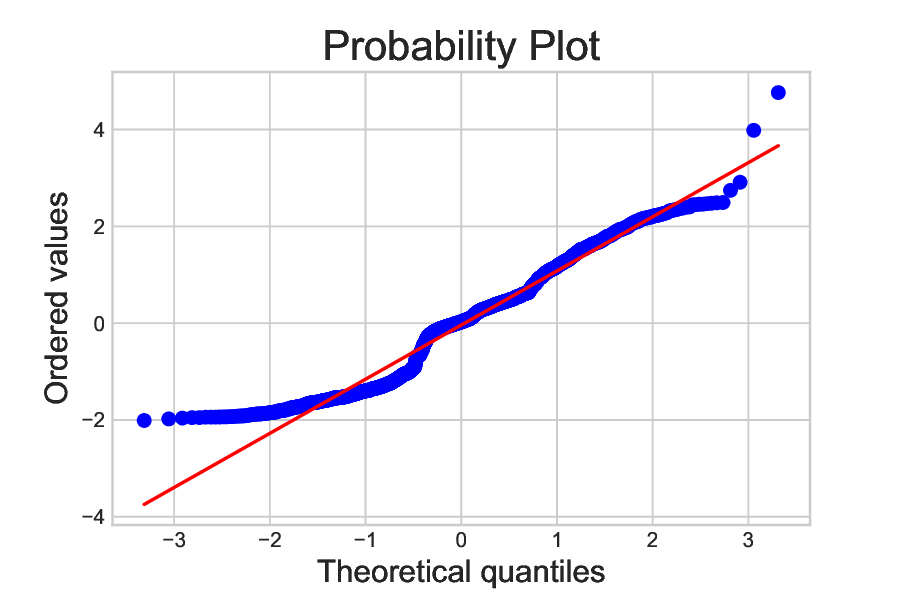}
	\caption{Normal probability plot of $H^{s}_{1}(N)$ on the second coordinate.}
	\label{fig:droite_henry_y}
   \end{minipage} 
\end{figure}

A visual tool to verify the normality of $H^{s}_{1}(N)$ for a particle is a
normal probability plot. Figures \ref{fig:droite_henry_x} and \ref%
{fig:droite_henry_y} displays the normal probability plot of $H^{s}_{1}(N)$
respectively for the $x$ axis and $y$ axis. For each axis, the normality is
confirmed: $H^{s}_{1}(N)$ fits well the theoretical quantiles. 

\begin{figure}[]
\centering
\includegraphics[width=0.6\textwidth]{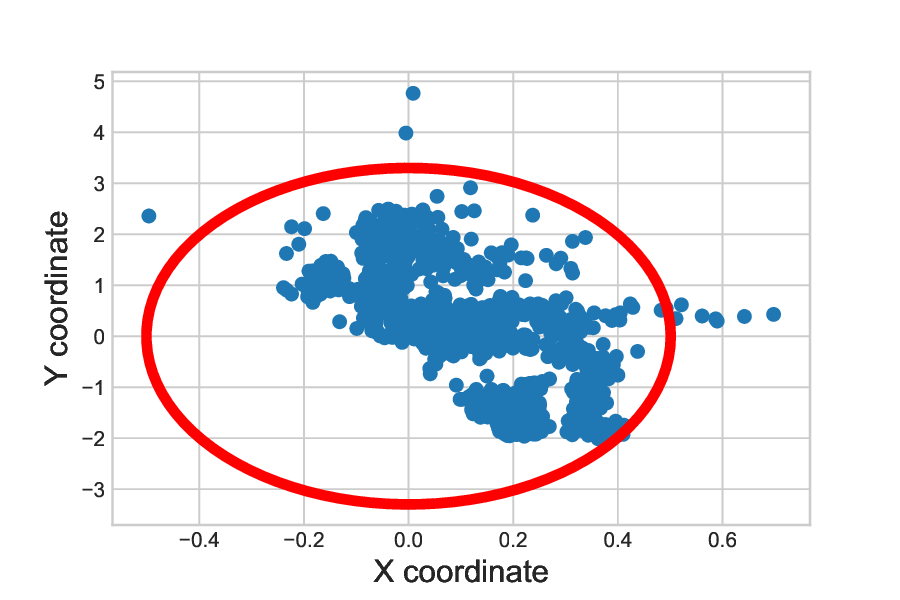}  
\caption{Trajectory of $H^{s}_{1}(N)$ for an oscillating particle in $\left[%
-10,10\right]^{2}$. The confidence ellipsoid at a level of 85\% is displayed
in red. Around 99\% of the trajectory of $H^{s}_{1}(N)$ is inside the
ellipse.}
\label{fig:ellipsoid}
\end{figure}

\begin{figure}[]
\centering
\includegraphics[width=\textwidth]{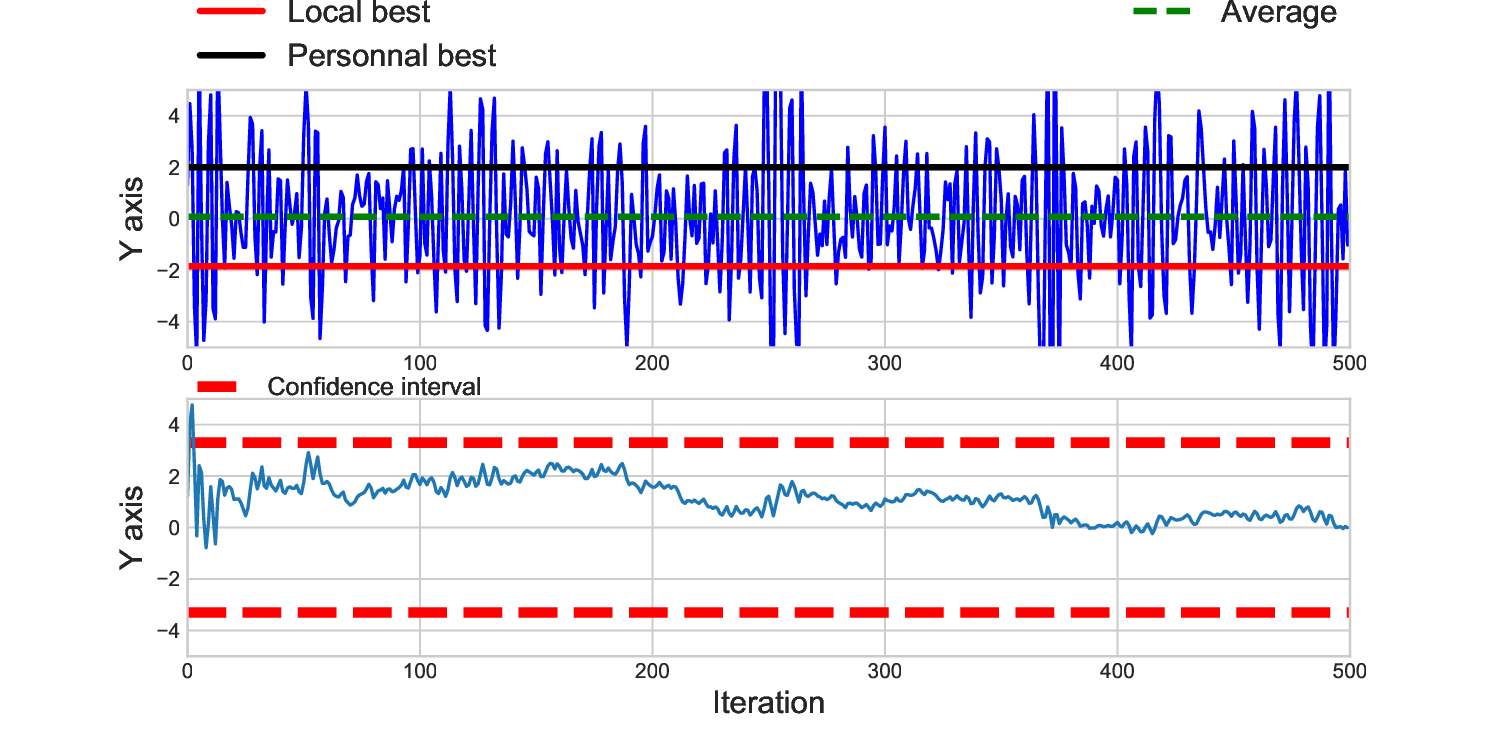}  
\caption{Top track: Trajectory of a oscillating particle on the y axis. The
particle is oscillating between 2 and -1.84. Bottom track: corresponding
trajectory of $H^{s}_{1}(N)$ on the y axis, the red dot are corresponding to
the 95\% confidence interval. The trajectory of $H^{s}_{1}(N)$ is well
bounded.}
\label{fig:trajectory_particle}
\end{figure}

To check the formula of the covariance matrix $\Gamma$, the confidence
ellipsoid is also a good indicator (see Figure \ref{fig:ellipsoid}). For a
single particle, $H^{s}_{1}(N)$ is not necessarily always inside the
confidence ellipsoid and does not respect the percentage of the defined
confidence level. Figure \ref{fig:trajectory_particle} shows the trajectory
of $x^{s}_n$ and $H^{s}_{1}(N)$ on the y axis, $H^{s}_{1}(N)$ remains
bounded. 

With 200 Monte-Carlo simulations of PSO (200 particles, and 2000
iterations), we select all the particles oscillating between $(3.58,-1.84)$
and $(3,2)$, and for each of them we compute $H^{s}_{1}(2000)$. Figure \ref%
{fig:density_mc_oscillating} displayed the density of $H^{s}_{1}(2000)$
using 1150 oscillating particles. Almost 95\% of the particles are inside
the confidence ellipsoid of level 95\% (represented in red). 
\begin{figure}[]
\centering
\includegraphics[width=0.75\textwidth]{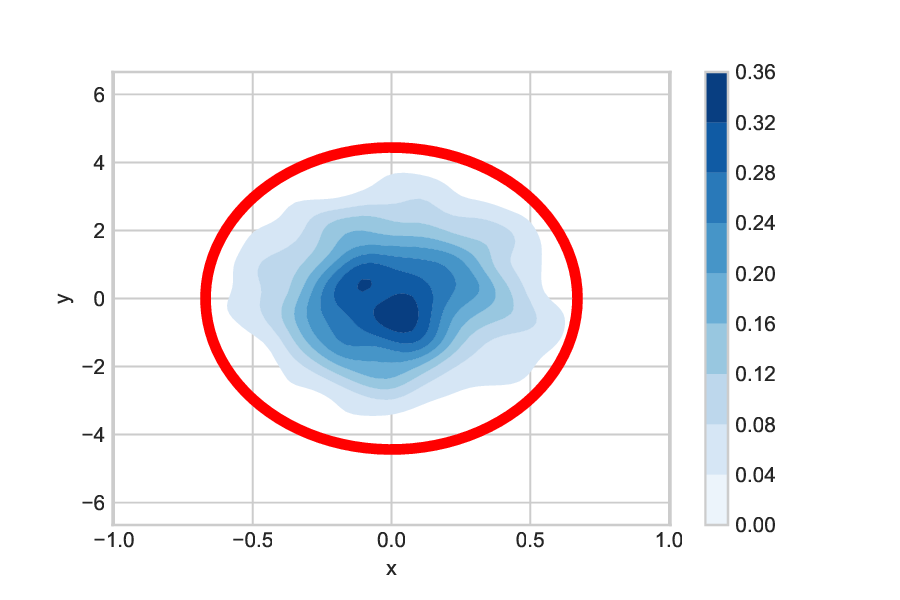}  
\caption{Density of $H_{1}(2000)$ with 1150 particles issued from
Monte-Carlo simulations. The red ellipse is the 95\% confidence ellipsoid. }
\label{fig:density_mc_oscillating}
\end{figure}
\FloatBarrier

\subsection{Non-oscillatory case}

\begin{figure}[]
\centering
\includegraphics[width=0.6\textwidth]{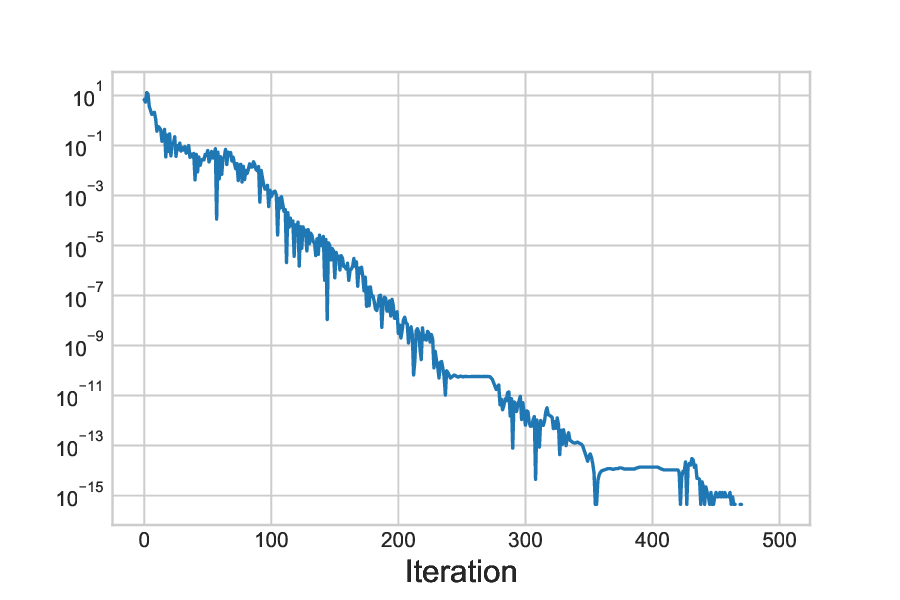}  
\caption{$\vert x_n-g_n\vert$ over 500 iterations in a logarithmic scale.
After 300 iterations, we reach the computer precision. We take advantage of
the linear behavior of $\log \left( \vert x_n-g_n \vert \right)$ on the 200
first iterations to perform a linear regression estimating $\protect\mu_x$.}
\label{fig:himmelblau_converging_linear}
\end{figure}

We study now the behaviors of non-oscillating particles on the Himmelblau's
function. We launch PSO with a population of 1000 particles and with 2000
iterations, $\omega=0.72984$ and $c=1.496172$. A ring topology was used to
ensure the presence of enough particle converging to each local optimum.  We
select particles converging to $(3,2)$, meaning that $p_n=g_n=p$ for a
sufficiently large $n$. For the weak convergence of the particle, we
consider: 
\[
H^{s}_{2}(N)=\frac{1}{\sqrt{N}}\left( \log\left\vert
x^{s}_{N}-g_{N}\right\vert -N \mu_{x}\right).
\]
First, it is easy to check the linear dependency of $\log\left\vert
x^{s}_{N}-g_{N}\right\vert$ with a single display of the trajectory. Figure %
\ref{fig:himmelblau_converging_linear} illustrates this phenomenon for a
single particle. We observed numerical issues when we reach the machine
precision, but a numerical approximation of $\mu_{x}$ can be performed
thanks to a linear regression.

\begin{figure}[]
\centering
\includegraphics[width=0.5\textwidth]{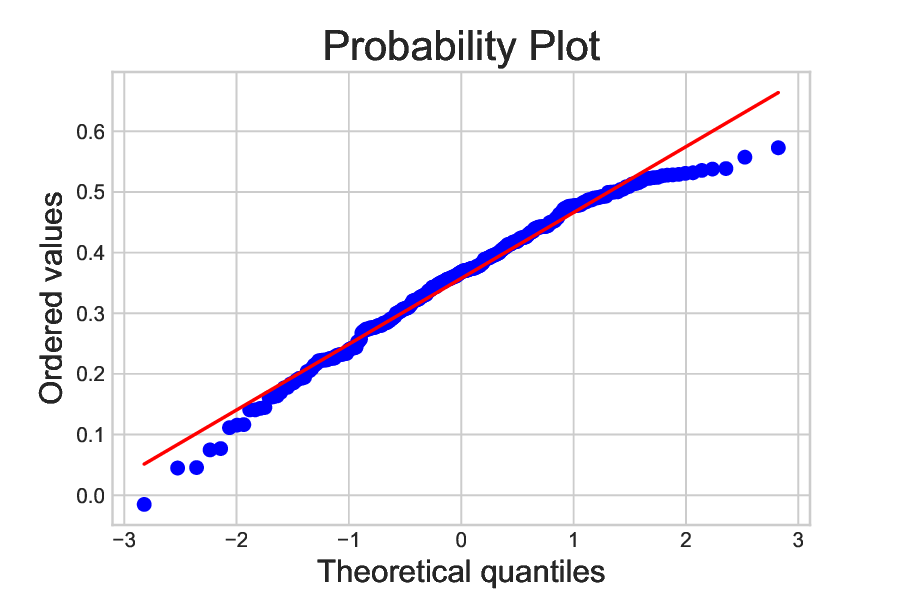}  
\caption{Normal probability plot of $H_{2}(N)$ on the first coordinate.}
\label{fig:henry_x_conv}
\end{figure}

\begin{figure}[]
\centering
\includegraphics[width=0.6\textwidth]{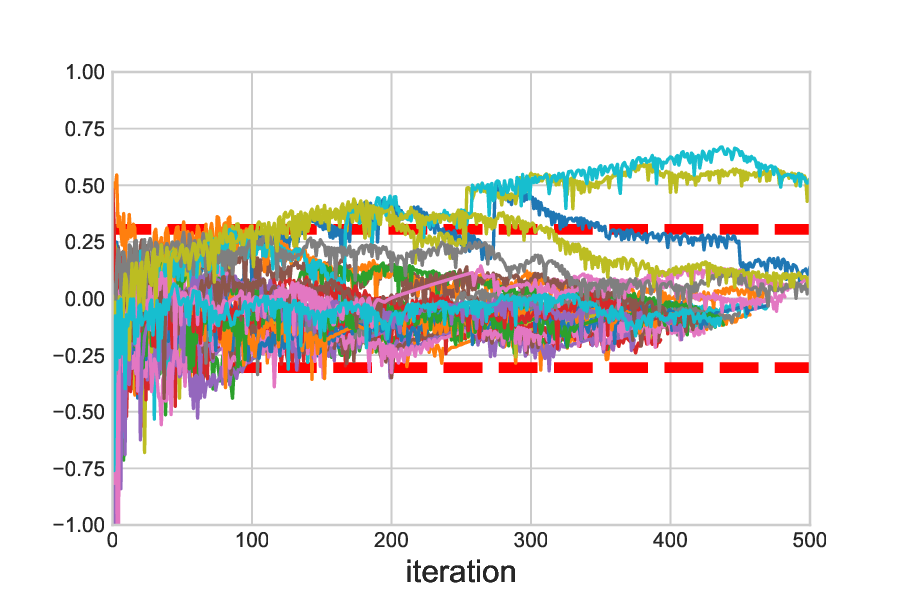}  
\caption{Trajectories of $H_{2}(N)$ on the first coordinate for five
particles. The red dot represents the 95\% confidence interval deduced from $%
\bar{\protect\sigma}_{x}$. Trajectories stop around the 400 iterations
(after an heating phase) due to numerical precision.}
\label{fig:multiple_H2}
\end{figure}

Using many converging particles, a Monte Carlo approximation of $\mu_x$ is
done. For the approximation of $\sigma_x$, a possibility is: 
\[
\bar{\sigma}_{x}^{2} =\mathrm{Var}_{\pi}\left( \log\left\vert
X_{0}\right\vert \right) +2\sum_{k=1}^{T}\mathrm{Cov}_{\pi}\left(
\log\left\vert X_{0}\right\vert ,\log\left\vert X_{k}\right\vert \right),
\]
where $T=20$. With near 240 converging particles to $(3,2)$, we found that
for the first coordinate: 
\begin{align*}
\bar{\mu}_{x} & = -0.032 , \\
\bar{\sigma}_{x} & =0.156.
\end{align*}
We verify the asymptotic normality of $H_{2}(N)$ with a normal probability
plot using the approximation of $\mu_x$. Figure \ref{fig:henry_x_conv}
displays the normal probability plot of $H_{2}(N)$ on the first coordinate,
the theoretical quantiles are well fitted by $H_{2}(N)$. 
Figure \ref{fig:multiple_H2} illustrates different trajectories of $H_{2}(N)$
on the first coordinate which are bounded by the 95\% confidence interval
deduced from $\bar{\sigma}_{x}$.

\subsection{Swarm at a fixed step}

To check the Theorem \ref{prop3}, we study: 
\[
H^{n}_{3}(S)=\frac{1}{\sqrt{S} \widehat{\sigma }_{n}}\sum_{i=1}^{S}\left(
x_n^i-g_n\right).
\]

\begin{figure}[]
\centering
\begin{minipage}[t]{.45\linewidth}
		\includegraphics[width=\textwidth]{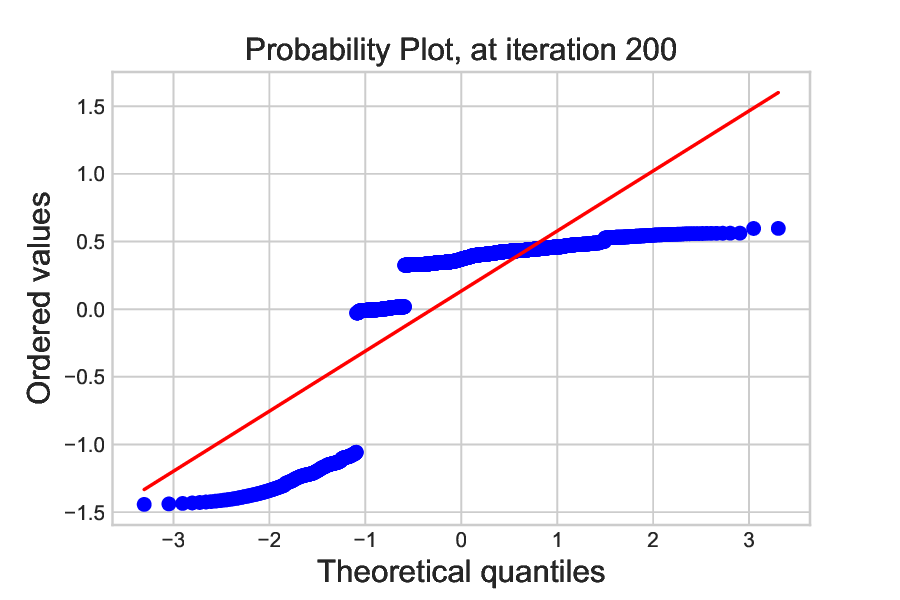}
	\caption{Normal probability plot of $H^{n=200}_{3}(S)$ at the 200th iteration on the first coordinate, using all the particles. We observe a discontinuity of the probability plot due to belated particles.}
	\label{fig:H3_all}
   \end{minipage}\hfill  
\begin{minipage}[t]{.45\linewidth}
		\includegraphics[width=\textwidth]{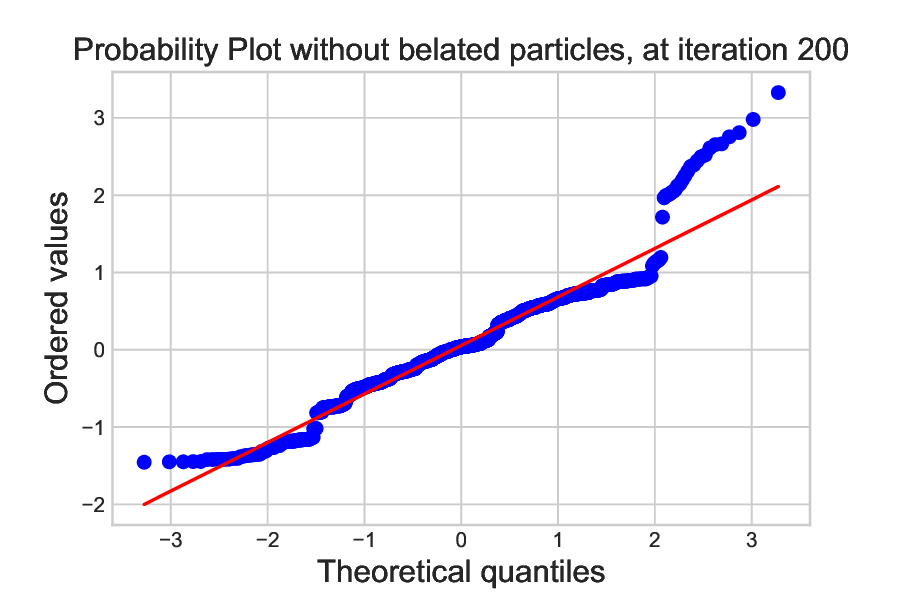}
	\caption{Normal probability plot of $H^{n=200}_{3}(N)$  at the 200th iteration on the first coordinate, computed without outlying or belated particles. }
	\label{fig:H3_sub}
   \end{minipage} 
\end{figure}

In practice, we encountered some difficulties to verify Theorem \ref{prop3}
because of the convergence rate of the particles. Indeed, when $%
p^s_n=g^s_n=p^s$, the particle $s$ converges exponentially to $g^s_n$ but
the spread of the rate of convergence is large. As a consequence, at a fixed
step of PSO, some particles could be considered as outliers because of a
lower rate of convergence. Because of these particles qualified as belated,
the asymptotic Gaussian behavior of $H^{n}_{3}(S)$ is not verified. A
solution is to filter the particles and remove the belated particles.
Figures \ref{fig:H3_all} and \ref{fig:H3_sub} illustrate this phenomenon for
the Himmelblau's function in 2D and with near 1500 converging particles to $%
(3,2)$ over 500 iterations. In Figure \ref{fig:H3_all}, we compute without
any filtering $H^{n}_{3}(S)$ and we notice that the Gaussian behavior is not
verified and some jumps appeared. The presence of these "jumps" is due to
belated particles which have a lower rate of convergence in comparison to
the swarm. When we remove these particles with a classical outliers
detection algorithm in Figure \ref{fig:H3_sub}, Theorem \ref{prop3} seems to
be verified.

\section{Conclusion}

Our main theoretical contribution revolves around the three CLTs and the
confidence regions derived from sampling either a particle's path or the
whole swarm. Practically the confidence set $\Lambda_{1-\alpha}$ localizes,
say $g$, with probability $1- \alpha$. The simulations carried out in Python language tend to foster the theory.

This work was initiated in order to solve a practical issue in connection
with oil industry and petrophysics. Yet in the previous section we confined
ourselves to simulated data for several reasons. It appears that our method
should be applied on real-data for a clear validation.

A second limitation of our work is the asymptotic set-up. The results are
stated for samples large enough. This may not have an obvious meaning for
people not familiar with statistics or probability. Practitioners usually
claim that the Central Limit Theorem machinery works well for a sample size
larger than thirty/forty whenever data are stable, say stationary (path with
no jumps...). As a consequence the first iterations of the algorithm should
be avoided to ensure a warming. When studying PSO the behavior of $p_n$ and $%
g_n$ turns out to be crucial too in order to ensure this stability, hence
the validity of the theoretical results. However the control of $p_n$ and $%
g_n$ is a difficult matter, which explains our current assumptions on
stagnation or ``almost stagnation". However, in order to fix the question of
the asymptotic versus non-asymptotic approach a work is on progress. We
expect finite-sample concentration inequalities : the non-asymptotic
counterparts of CLT.

In the oscillating framework our results
involve the parameter denoted $\theta$, the center of the segment $\left[p,g%
\right]$. Here we miss the target $g$ for instance. We claim that our method
may be adapted to tackle this issue but we will not elaborate on it. Briefly
speaking~: running two independent PSOs on the same domain with two distinct
pairs $\left(c_1,c_2\right)$ and $\left(c^{\prime}_1,c^{\prime}_2\right)$
will provide, under additional assumptions on the local and global minima of
the cost function, two CLTs involving : 
\[
\theta=\frac{c_1g+c_2p}{c_1+c_2} \quad \mathrm{and} \quad \theta^{\prime}=%
\frac{c^{\prime}_1g+c^{\prime}_2p}{c^{\prime}_1+c^{\prime}_2}. 
\]
A simple matrix inversion will then recover a joint CLT for the couple $%
\left(g,p\right)$ since, when $c_2 c^{\prime}_1 \neq c_1 c^{\prime}_2$ then :
\[
p=\frac{c c^{\prime}_1 \theta -c_1 c^{\prime}\theta^{\prime}}{c_2 c^{\prime}_1- c_1 c^{\prime}_2 }
\]
for instance.
\FloatBarrier

\section{Derivations of the results}

We start with some notations. First we recall that the sup-norm for square
matrices of size $d$ defined by: 
\[
\left\Vert M\right\Vert _{\infty}=\sup_{x\neq0}\frac{\left\Vert
Mx\right\Vert }{\left\Vert x\right\Vert }. 
\]

The tensor product notation is appropriate when dealing with special kind of
matrices for instance covariance matrices. Let $u$ and $v$ be two vectors of 
$\mathbb{R}^{d}$ then $u\otimes v=uv^{t}$ (where $v^{t}$ is the transpose of
vector $v$) stands for the rank one matrix defined for all vector $x$ by $%
\left( u\otimes v\right) \left( x\right) =\left\langle v,x\right\rangle u$.
Besides $\left\Vert u\otimes v\right\Vert _{\infty}=\left\Vert u\right\Vert
\left\Vert v\right\Vert $. The Hadamard product between vectors was
mentioned earlier. Its matrix version may be defined a similar way. Let $M$
and $S$ be two matrices with same size then $M\odot S$ is the matrix whose $%
\left( i,j\right) $ cell is $\left( M\odot S\right) _{i,j}=m_{i,j}s_{i,j}$.
We recall without proof the following computation rule mixing Hadamard and
tensor product. Let $\eta,\varepsilon,u$ and $v$ be four vectors in $\mathbb{%
R}^{d}$. Then: 
\begin{equation}
\left( \eta\odot u\right) \otimes\left( \varepsilon\odot v\right) =\left(
\eta\otimes\varepsilon\right) \odot\left( u\otimes v\right),  \label{Hadam}
\end{equation}
and the reader must be aware that the Hadamard product on the left-hand side
operates between vectors whereas on the right-hand side it operates on
matrices.

We will need $\mathcal{F}_{n}^{s}$ the filtration generated by the path of
particle number $s$ up to step $n$~: $\left\{x_{0}^{s},...,x_{n}^{s}\right\}$
and $\mathcal{F}_{n}^{\mathcal{S}}$ the filtration generated by the swarm up
to step $n$~: $\left\{ x_{0}^{s},...,x_{n}^{s}:s=1,...,S\right\} $ .

We will also denote later $g_{n}=\mathbb{E}\left(g_{n}\right)+\xi_{n}$ and $%
p_{n}=\mathbb{E}\left(p_{n}\right) +\nu_{n}$ the expectation-variance
decomposition of $g_{n}$ and $p_{n}$ where $\xi_{n}$ and $\nu_{n}$ are
centered random vectors and support all the variability of $g_{n}$ and $p_{n}
$ respectively.

\subsection{First case: oscillatory}

In this subsection we prove Theorem \ref{TH1}.
We start with two Lemmas and a Proposition who will be invoked later.

\begin{lemma}
\label{noise}Let $\varepsilon^{\left(i\right)}_{n}$ and $\eta^{\left(j\right)}_{n}$ be any coordinate of the random vectors $\varepsilon_{n}$ and $\eta_{n}$ appearing in \ref{base}. Clearly,  $\varepsilon^{\left(i\right)}_{n}$ and $\eta^{\left(j\right)}_{n}$  are not
independent but not correlated and follow the same type of distribution. Besides: 
\begin{align*}
\mathbb{E}{\varepsilon^{\left(i\right)}}^{2}  &  =\mathbb{E}{\eta^{\left(j\right)}}^{2}=1/6,\quad\mathbb{E}
\varepsilon^{\left(i\right)}\eta^{\left(j\right)}=0,\\
\mathbb{E}{\varepsilon^{\left(i\right)}}^{3}\eta^{\left(j\right)} &  =\mathbb{E}\varepsilon^{\left(i\right)}{\eta^{\left(j\right)}}^{3}=0.
\end{align*}
\end{lemma}
The proof is very simple hence omitted.

The next Lemma is also straigthforward but will be frequetly invoked later. 
\begin{lemma}
\label{LLN}Let $E_{n}$ be a sequence of i.i.d centered random matrices with
finite moment of order 4, let $u_{n}$ and $v_{n}$ two sequence of random
vectors almost surely bounded and such that $\left(  u_{n},v_{n}\right)  $ is
for all $n$ independent from $E_{n}$ then for the $\left\Vert \cdot
\right\Vert _{\infty}$ norm:
\[
\frac{1}{N}\sum_{n=1}^{N}E_{n}\odot\left(  u_{n}\otimes v_{n}\right)
\rightarrow_{\mathbb{P}}0.
\]
\end{lemma}
The proof is a simple application of Cauchy-Schwartz inequality.

\noindent \textbf{Proof of Theorem \ref{TH1}:}

\textbf{First step :} \textit{We show that the problem of convergence in
distribution for }$\frac{1}{N}\sum_{n=1}^{N}\left(  x_{n}-\theta\right)
$\textit{ may be shifted to proving a CLT for a martingale difference array.}

Starting from the initial PSO equation :%
\[
x_{n+1}-x_{n}=\omega\left(  x_{n}-x_{n-1}\right)  +c_{1}r_{1,n}\odot\left(
p_{n}-x_{n}\right)  +c_{2}r_{2,n}\odot\left(  g_{n}-x_{n}\right)
\]
it is simple to derive :%
\[
x_{n+1}-x_{n}=\omega\left(  x_{n}-x_{n-1}\right)  +c_{1}r_{1,n}\odot\left(
p-x_{n}\right)  +c_{2}r_{2,n}\odot\left(  g-x_{n}\right)  +R_{n}%
\]
with
\[
R_{n}=c_{1}r_{1,n}\odot\left(  p_{n}-p\right)  +c_{2}r_{2,n}\odot\left(
g_{n}-g\right)  .
\]
Denote $z_{n}=\left(  x_{n}-\theta\right)  $ then after some additional
calculations we get :%
\begin{equation} \label{decomp-base}
z_{n+1}-\left(  1+\omega-c\right)  z_{n}+\omega z_{n-1}=c_{1}\varepsilon
_{1,n}\odot\left(  \frac{c_{2}}{c_{1}+c_{2}}\left(  p-g\right)  -z_{n}\right)
-c_{2}\varepsilon_{2,n}\odot\left(  \frac{c_{1}}{c_{1}+c_{2}}\left(
p-g\right)  +z_{n}\right)  +R_{n}%
\end{equation}
with, for $i=1,2:$ $r_{i,n}=\mathbb{E}r_{i,n}+\varepsilon_{i,n}%
=1/2+\varepsilon_{i,n}$ and $\varepsilon_{i,n}$ is the centered version of
$r_{i,n}$.

Summing the above equation from $n=1$ to $n=N-1$ :%
\begin{equation}
\sum_{n=1}^{N}z_{n}=\frac{1}{c}\sum_{n=1}^{N}u_{n}+\mathcal{R}%
_{N}\label{decomp-zn-un}%
\end{equation}
with :%
\begin{align*}
u_{n} &  =c_{1}\varepsilon_{1,n}\odot\left(  \frac{c_{2}}{c_{1}+c_{2}}\left(
p-g\right)  -z_{n}\right)  -c_{2}\varepsilon_{2,n}\odot\left(  \frac{c_{1}%
}{c_{1}+c_{2}}\left(  p-g\right)  +z_{n}\right),  \\
\mathcal{R}_{N} &  =\frac{1}{c}\left[  \sum_{n=1}^{N}R_{n}+z_{1}+\left(
c-1\right)  z_{N}+\omega\left(  z_{N-1}-z_{0}\right)  \right]
\end{align*}
Notice that $u_{n}$ just above is a martingale difference with respect to the nested
filtration $\mathcal{F}_{n}$ : $\mathbb{E}\left(  u_{n}|\mathcal{F}%
_{n-1}\right)  =0$ because $\varepsilon_{1,n}$ and $\varepsilon_{2,n}$ are
centered and independent from $z_{n}$.

In view of conditions $\mathbf{A}_{1}$ and $\mathbf{A}_{2}$  :%
\[
\frac{\mathcal{R}_{N}}{\sqrt{N}}\rightarrow_{\mathbb{P}}0,
\]

then from (\ref{decomp-zn-un}) $\frac{1}{\sqrt{N}}\sum_{n=1}^{N}z_{n}$
converges in distribution if and only if $\frac{1}{\sqrt{N}}\sum_{n=1}%
^{N}u_{n}$ does.

\textbf{Second step }: \textit{it is shown below that weak convergence for
}$\frac{1}{\sqrt{N}}\sum_{n=1}^{N}u_{n}$\textit{ holds under the assumptions
of the Theorem.}

We aim at proving a Levy-Lindeberg version
of the CLT for the series of $u_{n}$ in two steps (Theorem 2.1.9 p. 46 and
its corollary 2.1.10 in \cite{Duflo_1997}): first checking the Lyapunov condition holds (hence the
Lindeberg's uniform integrability that ensures uniform tightness of the
sequence) then ensuring convergence in probability of the conditional covariance structure of $\left( 1/\sqrt{N}\right)
\sum_{n=1}^{N}u_{n}$.

Here the Lyapunov condition holds trivially because we are faced with bounded martingale difference sequences. Indeed by
assumption $\mathbf{A}_{1}$ and since the random variables $\varepsilon_{1,n}$ and
$\varepsilon_{2,n}$ are almost surely bounded $u_{n}$ is itself almost surely
bounded.

We turn to the conditional covariance sequence of $u_n$ :%
\[
\Gamma_{N}=\frac{1}{N}\sum_{n=1}^{N}\mathbb{E}\left(  u_{n}\otimes
u_{n}|\mathcal{F}_{n-1}\right).
\]

Starting from the definition of $u_n$ just below (\ref{decomp-zn-un}) leads after simple calculations to :%
\begin{align*}
\mathbb{E}\left(  u_{n}\otimes u_{n}|\mathcal{F}_{n-1}\right)  &
=\frac{c_{1}^{2}c_{2}^{2}}{\left(  c_{1}+c_{2}\right)  ^{2}}\frac{1}%
{6}\mathbf{I}\odot\left[  \left(  p-g\right)  \otimes\left(  p-g\right)
\right]  +\frac{\left(  c_{1}^{2}+c_{2}^{2}\right)  }{12}\mathbf{I}%
\odot\left(  z_{n}\otimes z_{n}\right)  \\
&  -2\frac{c_{1}-c_{2}}{12}\frac{c_{1}c_{2}}{c_{1}+c_{2}}\mathbf{I}%
\odot\left[  \left(  p-g\right)  \otimes z_{n}\right]
\end{align*}

Notice that, due to the operation $\mathbf{I}\odot\left(.\right)$  $\mathbb{E}\left( u_{n}\otimes u_{n}|\mathcal{F}_{n-1}\right) $
is a diagonal matrix with :%
\begin{align*}
\mathbb{E}\left(  u_{n}\otimes u_{n}|\mathcal{F}_{n-1}\right)   &
=diag(\tau_{i,n})_{1\leq i\leq d}\\
\tau_{i,n} &  =\frac{c_{1}^{2}c_{2}^{2}}{6\left(  c_{1}+c_{2}\right)  ^{2}%
}\left(  p_{i}-g_{i}\right)  ^{2}+\frac{\left(  c_{1}^{2}+c_{2}^{2}\right)
}{12}z_{i,n}^{2}-\frac{c_{1}-c_{2}}{6}\frac{c_{1}c_{2}}{c_{1}+c_{2}}\left[
\left(  p_{i}-g_{i}\right)  z_{in}\right]
\end{align*}
Averaging the previous results in $n$ we get :%
\begin{align*}
\Gamma_{N} &  =\frac{c_{1}^{2}c_{2}^{2}}{\left(  c_{1}+c_{2}\right)  ^{2}%
}\frac{1}{6}\mathbf{I}\odot\left[  \left(  p-g\right)  \otimes\left(
p-g\right)  \right]  +\frac{\left(  c_{1}^{2}+c_{2}^{2}\right)  }%
{12}\mathbf{I}\odot\left(  \frac{1}{N}\sum_{n=1}^{N}z_{n}\otimes z_{n}\right)
\\
&  -\frac{c_{1}-c_{2}}{6}\frac{c_{1}c_{2}}{c_{1}+c_{2}}\mathbf{I}\odot\left[
\left(  p-g\right)  \otimes\left(  \frac{1}{N}\sum_{n=1}^{N}z_{n}\right)
\right]
\end{align*}

The last term is negligible : indeed Lemma \ref{LLN} shows that $\frac{1}{N}\sum_{n=1}^{N}z_{n}%
\rightarrow_{\mathbb{P}}0$. Our only task consists in
studying convergence for $\mathbf{I}\otimes \left(\frac{1}{N}\sum_{n=1}^{N}z_{n}\otimes z_{n}\right)$. More
precisely we just need to prove convergence for $\left(  1/N\right)
\sum_{n=1}^{N}z_{i,n}^{2}$ for all $i\in\left\{ 1,...,d \right\}  $.

\textbf{Third (and last) step :} \textit{Proof of convergence of }$\frac{1}%
{N}\sum_{n=1}^{N}z_{in}^{2}$ for all $i$.

Let us restart from (\ref{decomp-base}) :%
\[
z_{n+1}-\left(  1+\omega-c\right)  z_{n}+\omega z_{n-1}=u_{n}+R_{n}%
\]
Clearly $R_{n}$ is a negligible term due to assumption $\mathbf{A}_2$. All the authors' computations were effectively carried out with $R_n$ but in order to alleviate this proof we will proceed here as if :%
\[
z_{i,n+1}=\left(  1+\omega-c\right)  z_{i,n}-\omega z_{i,n-1}+u_{i,n}%
\]
From now on we drop the dimension index $i$ in the previous equation :%
\begin{align*}
z_{i,n+1}^{2}  &  =\left(  1+\omega-c\right)  ^{2}z_{i,n}^{2}+\omega^{2}%
z_{i,n-1}^{2}+u_{i,n}^{2}\\
&  +2\left[  \left(  1+\omega-c\right)  z_{i,n}-\omega z_{i,n-1}\right]
u_{i,n}-2\omega\left(  1+\omega-c\right)  \left[  z_{i,n}z_{i,n-1}\right]
\end{align*}
The following equation is derived from the preceding one after summing and
considering the $i^{th}$ component of each vector :%
\begin{equation}\label{z2}
\begin{aligned}
\left(  1-\left(  1+\omega-c\right)  ^{2}-\omega^{2}\right)  \sum_{n=1}%
^{N}z_{i,n}^{2}  = &\sum_{n=1}^{N}u_{i,n}^{2}-2\omega\left(  1+\omega
-c\right)  \sum_{n=1}^{N}z_{i,n}z_{i,n-1} \\
 & +2\sum_{n=1}^{N}\left[  \left(  1+\omega-c\right)  z_{i,n}-\omega
z_{i,n-1}\right]  u_{i,n}+z_{i,1}^{2}-z_{i,N+1}^{2}+\omega^{2}\left(
z_{i,0}^{2}-z_{i,N}^{2}\right)
\end{aligned}
\end{equation}

It is easily seen that Lemma \ref{LLN} may be applied to the series with term $\left[  \left(  1+\omega-c\right)
z_{i,n}-\omega z_{i,n-1}\right]  u_{i,n}$ and its average tends to zero.

From now on we remove all indices $i$ for components of the vectors in order
to alleviate notations.

Our next task consists in studying $\sum_{n=1}^{N}z_{n}z_{n-1}$. To that
purpose we turn back again to equation (\ref{decomp-base}) and multiply by $z_{n}$ the whole
equation and sum over $n$ :%
\[
\sum_{n=1}^{N}z_{n+1}z_{n}-\left(  1+\omega-c\right)  \sum_{n=1}^{N}z_{n}%
^{2}+\omega\sum_{n=1}^{N}z_{n}z_{n-1}=\sum_{n=1}^{N}z_{n}u_{n}
\]
hence%
\[
\left(  1+\omega\right)  \sum_{n=1}^{N}z_{n}z_{n-1}=\left(  1+\omega-c\right)
\sum_{n=1}^{N}z_{n}^{2}+\sum_{n=1}^{N}z_{n}u_{n}+z_{1}z_{0}-z_{N+1}z_{N}%
\]
The same arguments as above show that $\left(  1/N\right)  \sum_{n=1}^{N}%
z_{n}u_{n}$ tends to $0$ in probability by the weak law of large numbers. This proves that when $N \to +\infty$ :%
\[
\frac{1}{N}\left(  \sum_{n=1}^{N}z_{n}z_{n-1}-\frac{1+\omega-c}{1+\omega
}\sum_{n=1}^{N}z_{n}^{2}\right)  \rightarrow_{\mathbb{P}}0
\]
Plugging this result in (\ref{z2}) and taking into account the remarks above :%
\begin{equation}\label{jdg}
\lim_{N \to +\infty}\frac{1}{N}\left\{  \left[  \left(  1-\left(  1+\omega-c\right)
^{2}-\omega^{2}\right)  +2\omega\frac{\left(  1+\omega-c\right)  ^{2}%
}{1+\omega}\right]  \sum_{n=1}^{N}z_{n}^{2}-\sum_{n=1}^{N}u_{n}^{2}\right\}
=0
\end{equation}
Our last task consists in finding $\lim_{N \to +\infty}\left(  1/N\right)  \sum_{n=1}^{N}%
u_{n}^{2}$. Calculations are closed to those carried out at the beginning of the second step but this time without conditionning with respect to $\mathcal{F}_n$. Remind that the index $i$ was removed for the sake of clarity :
\begin{equation}\label{jdg2}
u_{n}^{2}  =\left[c_{1}\varepsilon_{1,n} \left(  \frac{c_{2}}{c_{1}+c_{2}}\left(
p-g\right)  -z_{n}\right)  -c_{2}\varepsilon_{2,n} \left(  \frac{c_1}{c_1 + c_2}\left(  p-g\right)  +z_{n}\right) \right]  ^{2} \\
\end{equation}

Developping the square into bracket and taking into account that the cross product will tend to zero after averaging namely :
\[
\lim_{N \to +\infty }\frac{1}{N} \sum_{n=1}^{N} \varepsilon_{1,n}\varepsilon_{2,n} \left(  \frac{c_{2}}{c_{1}+c_{2}}\left(
p-g\right)  -z_{n}\right)\left(  \frac{c_1}{c_1 + c_2}\left(  p-g\right)  +z_{n}\right) =0 
\]

in probability, we can focus only on : 
\[
v_{1,N}=\frac{c_{1}^2}{N} \sum_{n=1}^{N}  \varepsilon_{1,n}^2 \left(  \frac{c_{2}}{c_{1}+c_{2}}\left(
p-g\right)  -z_{n}\right)^2 \quad \mathrm{and} \quad  v_{2,N}=\quad \frac{c_{2}^2}{N} \sum_{n=1}^{N}  \varepsilon_{2,n}^2 \left(  \frac{c_{1}}{c_{1}+c_{2}}\left(p-g\right) +z_{n}\right)^2 
\]

Both $v_{1,N}$ and $v_{2,N}$ may be managed the same way. We develop only for $v_{1,N}$ :
\[
v_{1,N}=\frac{c_{1}^2}{N} \sum_{n=1}^{N}  \left(\varepsilon_{1,n}^2 - \mathbb{E}\varepsilon_{1,n}^2\right)\left(  \frac{c_{2}}{c_{1}+c_{2}}\left(
p-g\right)  -z_{n}\right)^2 + \frac{c_{1}^2}{N} \sum_{n=1}^{N}  \mathbb{E}\varepsilon_{1,n}^2 \left(  \frac{c_{2}}{c_{1}+c_{2}}\left(
p-g\right)  -z_{n}\right)^2
\]
The first series above tends in probability to 0 again by the weak law of large numbers for martingale difference arrays already mentioned several times earlier in the proof. The second series depends only on the $z_i$'s.
We finally get :
\[
\lim_{N \to +\infty} v_{1,N}=\lim_{N \to +\infty} \frac{c_{1}^2}{12 N}\sum_{n=1}^{N}  \left(  \frac{c_{2}}{c_{1}+c_{2}}\left(
p-g\right)  -z_{n}\right)^2
\]

Similarly :
\[
\lim_{N \to +\infty} v_{2,N}=\lim_{N \to +\infty} \frac{c_{2}^2}{12 N}\sum_{n=1}^{N}  \left(  \frac{c_{1}}{c_{1}+c_{2}}\left(
p-g\right)  +z_{n}\right)^2
\]

Putting the two last equations with (\ref{jdg}) and (\ref{jdg2}) we conclude that :

%\\
%& = \left[c_{1}\varepsilon_{1,n} \left(  \frac{c_{2}}{c_{1}+c_{2}}\left(
%p-g\right)  -z_{n}\right)\right]  ^{2} + \left[c_{2}\varepsilon_{2,n} \left(  \frac{c_{1}%
%}{c_{1}+c_{2}}\left(  p-g\right)  +z_{n}\right) \right]  ^{2}
 %\\
%& -2 c_{1} c_{2}  \varepsilon_{1,n} \varepsilon_{2,n}\left(  \frac{c_{2}}{c_{1}+c_{2}}\left(
%p-g\right)  -z_{n}\right) \left(  \frac{c_{1} }{c_{1}+c_{2}}\left(  p-g\right)  +z_{n}\right) 

%u_{n}^{2}  &  =\left[  \frac{c_{1}c_{2}}{c_{1}+c_{2}}\left(  \varepsilon
%_{1,n}-\varepsilon_{2,n}\right)  \left(  p-g\right)  -\left(  c_{1}%
%\varepsilon_{1,n}+c_{2}\varepsilon_{2,n}\right)  z_{n}\right]  ^{2}\\
%&  =\left[  \frac{c_{1}c_{2}}{c_{1}+c_{2}}\left(  \varepsilon_{1,n}%
%-\varepsilon_{2,n}\right)  \left(  p-g\right)  \right]  ^{2}+\left(
%c_{1}\varepsilon_{1,n}+c_{2}\varepsilon_{2,n}\right)  ^{2}z_{n}^{2}%
%-2\frac{c_{1}c_{2}}{c_{1}+c_{2}}\left(  \varepsilon_{1,n}-\varepsilon
%_{2,n}\right)  \left(  c_{1}\varepsilon_{1,n}+c_{2}\varepsilon_{2,n}\right)
%\left(  p-g\right)  z_{n}%

\[
\lim_{N \to +\infty}\left\{  \left(  1-\left(  1+\omega-c\right)  ^{2}-\omega^{2}\right)
+2\omega\frac{\left(  1+\omega-c\right)  ^{2}}{1+\omega}-\frac{c_{1}^{2}%
+c_{2}^{2}}{12}\right\}  \frac{1}{N}\sum_{n=1}^{N}z_{i,n}^{2}=\frac{1}%
{6}\left(  \frac{c_{1}c_{2}}{c_{1}+c_{2}}\right)  ^{2}\left(  p-g\right)  ^{2}%
\]%

However the previous equation is valid only if the big term between brackets in the right hand side is strictly positive. This term depends only on $c_1$, $c_2$ and $\omega$. This compatibility condition ensures just that a non-null limit for $\sum_{n=1}^{N}z_{i,n}^{2}$ exists. 

Besides :

\[ 1-\left(  1+\omega-c\right)  ^{2}-\omega^{2}+2\omega\frac{\left(
1+\omega-c\right) ^{2}}{1+\omega}-\frac{c_{1}^{2}+c_{2}^{2}}{12}
=  2c\frac{1-\omega}{1+\omega}\left(  1+\omega-\frac{c}{2}\right)
-\frac{c_{1}^{2}+c_{2}^{2}}{12}
\]
and the positivity condition matches assumption $\mathbf{A}_3$ within the Theorem.

This yields :
\[
\lim_{N \to +\infty}\frac{1}{N}\sum_{n=1}^{N}z_{i,n}^{2}=\frac{1}{6}\left(  \frac{c_{1}c_{2}%
}{c_{1}+c_{2}}\right)  ^{2}\left\{  2c\frac{1-\omega}{1+\omega}\left(
1+\omega-\frac{c}{2}\right)  -\frac{c_{1}^{2}+c_{2}^{2}}{12}\right\}
^{-1}\left(  p-g\right)  ^{2}%
\]

hence from (\ref{}) and with simple algebra :

\[
\Gamma_N \rightarrow_{\mathbb{P}} \mathfrak{K} \mathrm{diag}\left(  p-g\right)^{\odot2}
\]
with
\[ \mathfrak{K}=
\frac{1}{24} \left(\frac{c_{1}c_{2}}{c}\right)^2
\left[
2c\frac{1-\omega}{1+\omega}\left(1+\omega-\frac{c}{2}\right)  
\left\{  2c\frac{1-\omega}{1+\omega}\left(
1+\omega-\frac{c}{2}\right)  -\frac{c_{1}^{2}+c_{2}^{2}}{12}\right\}
\right]^{-1}
\]

All that was done until this point ensures that $\frac{1}{N} \sum_{n=1}^N u_n \hookrightarrow \mathbf{Z}$ where $\mathbf{Z}$ is a centered Gaussian random vector in $\mathbb{R}^d$ with covariance operator $\mathfrak{K} \mathrm{diag}\left(  p-g\right)^{\odot2}$.

At last (\ref{decomp-zn-un}) leads to the desired result : $\frac{1}{N} \sum_{n=1}^N z_n \hookrightarrow \mathbf{Z}/c$.
\endproof

\subsection{Second case: non-oscillatory and stagnant}

We start from (\ref{chain:X}) and set $X_{n}=y_{n}/y_{n-1}$. Then we get : 
\[
X_{n+1}=1+\omega-c-\frac{\omega}{X_{n}}+c\varepsilon_{n}, 
\]
where $\varepsilon_{n}=r_{1,n}+r_{2,n}-1$ has a ``witch hat" distribution
(convolution of two uniform distributions) with support $\left[ -1,+1\right] 
$. We focus now on the above homogeneous Markov chain $X_{n}$ and we aim at
proving that a CLT holds for $h\left( X_{n}\right) =\log\left\vert
X_{n}\right\vert $ namely that for some $\mu$ and $\sigma^{2}$ : 
\[
\sqrt{N}\left[ \frac{1}{N}\sum_{n=1}^{N}\log\left\vert X_{n}\right\vert -\mu%
\right] \hookrightarrow\mathcal{N}\left( 0,\sigma^{2}\right), 
\]
which will yield: 
\[
\sqrt{N}\left[ \frac{1}{N}\log\left\vert y_{N}\right\vert -\mu\right]
\hookrightarrow\mathcal{N}\left( 0,\sigma^{2}\right). 
\]

We aim at applying Theorem 1 p. 302 in \cite{jones2004markov} (see also \cite%
{meyn2012markov}, section 17.5 for similar results). We need to check three
points : (i) $X_n$ is Harris ergodic, (ii) the existence of a small set $%
\mathcal{C}$ and (iii) of a function $g$ with a drift condition (see \cite%
{meyn2012markov}) such that (\ref{maj_log}) below holds. Following \cite%
{Hairer2011}, section 3.2 it turns out that proving (ii) and (iii) above
with a drift function $g$ but without the minoration bound (\ref{maj_log})
is sufficient to ensure (i).

\bigskip

Denote $P\left( t,x\right) $ the transition kernel of $X_n$. It is plain
that $P\left( t,x\right) $ coincides with the density of the uniform
distribution on the set : 
\[
\mathcal{E}_{x}=\left[ 1+\omega -2c-\frac{\omega }{x },1+\omega -\frac{%
\omega }{x}\right]. 
\]
Theorem \ref{theo:equal} is a consequence of the two Lemmas below coupled
with the above-mentioned Theorem 1 p.302 in \cite{jones2004markov}.

\begin{lemma}
\label{smallset}Take $M_{\tau}=\omega/\left( c-\tau\right) $ with any $
0<\tau<c$ then the set $\mathcal{C}=\left( -\infty,-M_{\tau}\right] \cup
\left[ M_{\tau},+\infty\right) $ is a small set for the transition kernel of 
$X_{n}$.
\end{lemma}

\noindent\textbf{Proof :}

\noindent We have to show that for all $x\in \mathcal{C}$ and Borel set $A$
in $\mathbb{\ R}$: 
\[
P\left( A,x\right) \geq \varepsilon Q\left( A\right) ,
\]
where $\varepsilon >0$ and $Q$ is a probability distribution. The main
problem here comes from the compact support of $P\left( t,x\right) $. Take $%
x $ such that $\left\vert x\right\vert \geq M$ then: 
\[
1+\omega -c-\frac{\omega }{M}+c\varepsilon _{n}\leq 1+\omega -c-\frac{\omega 
}{x}+c\varepsilon _{n}\leq 1+\omega -c+\frac{\omega }{M}+c\varepsilon _{n},
\]
where $\varepsilon _{n}$ has compact support $\left[ -1,+1\right] $. It is
simple to see that with $M=M_{\tau }=\omega /\left( c-\tau \right) $ the
above bound becomes: 
\[
1+\omega -2c+\tau +c\varepsilon _{n}\leq 1+\omega -c-\frac{\omega }{x}
+c\varepsilon _{n}\leq 1+\omega +\tau +c\varepsilon _{n}.
\]
The intersection of the supports of $1+\omega -2c+\tau +c\varepsilon _{n}$
and $1+\omega+\tau +c\varepsilon _{n}$ is the set $\left[ 1+\omega -c-\tau
,1+\omega -c+\tau \right] $ whatever the value of $x$ in $\mathcal{C} $. The
probability measure $Q$ mentioned above may be chosen as the uniform
distribution with support $\left[ 1+\omega -c-\tau ,1+\omega -c+\tau \right] 
$.\bigskip

Now we turn to the drift condition. Our task consists in constructing a
function \newline
$g:\mathbb{R}\rightarrow\left[ 1,+\infty\right[ $ such that for all $x$: 
\begin{equation}
\int_{\mathbb{R}}g\left( t\right) P\left( t,x\right) dt\leq\rho _{1}g\left(
x\right) +\rho_{2}1_{x\in\mathcal{C}} ,  \label{drift_condition}
\end{equation}
where $0<\rho_{1}<1$ and $\rho_{2}\geq0.$ Besides, in order to get a CLT on $%
\log\left\vert X_{n}\right\vert $ we must further ensure that for all $x$: 
\begin{equation}
\log^{2}\left\vert x\right\vert\leq g\left( x\right).  \label{maj_log}
\end{equation}
Note however that, if (\ref{drift_condition}) holds for $g$ but (\ref%
{maj_log}) fails, then both conditions will hold for updated function $%
g^{\ast}=\eta g$ with constant $\eta>1$ and $\rho_{2}^{\prime}=\eta\rho_{2}$
such that (\ref{maj_log}) holds.

The next Lemma constructs the function $g$ mentioned above. 
\begin{lemma}
\label{drift}Take for $g$ the even function defined by $g\left( x\right)
=C_{1}/\sqrt{\left\vert x\right\vert }$ for $\left\vert x\right\vert \leq
M_{\tau}$ and $g\left( x\right) =C_{2}\left( \log\left\vert x\right\vert
\right) ^{2}$ for $\left\vert x\right\vert >M_{\tau}$. Assume that $\mathbf{B}_{1}$ holds. \\
Then it is always possible to choose three constants $\tau,$ $C_{1}$ and $
C_{2}$ such that (\ref{drift_condition}) holds for a specific choice of $
\rho_{1}$ and $\rho_{2}$.
\end{lemma}

The proof of Lemma \ref{drift} is postponed to the Appendix.

%\begin{acknowledgements}
%\end{acknowledgements}


%\bibliographystyle{plain}
\bibliography{pso_chapter_bib}

\newpage

\appendix

\appendixpage

This appendix provides a mathematical derivation of Lemma \ref{drift}.

\section{Proof of Lemma \protect\ref{drift}}

The proof of the Lemma just consists in an explicit construction of the
above-mentioned $\tau$, $C_{1}$, and $C_{2}$. This construction is detailed
for the sake of completeness.

\textbf{At this point and in order to simplify the computations below we
will assume that the distribution of }$\varepsilon_{n}$\textbf{\ is uniform
on }$\left[ -1,+1\right] $ instead of the convolution of two $\mathcal{U} _{ %
\left[ -1/2,1/2\right] }$ distributions.

Set $\lambda =1+\omega -c$, assume that $\lambda >0$ (the case $\lambda <0$
follows the same lines) and notice that: 
\[
\int_{\mathbb{R}}g\left( t\right) P\left( t,x\right) dt=\frac{1}{2c}\int_{ 
\mathbb{\lambda -}\left( \omega /x\right) -c}^{\mathbb{\lambda -}\left(
\omega /x\right) +c}g\left( s\right) ds=\frac{1}{2c}\int_{\left( \omega
/x\right) -\lambda -c}^{\left( \omega /x\right) -\lambda +c}g\left( s\right)
ds,
\]
the last inequality stemming from parity of $g$. We should consider two
cases $x>0$ and $x\leq 0$.

The proof takes 2 parts ($x>0$ and $x<0$ respectively). Both are given again
for completeness and because the problem is not symmetric. Each part is
split in three steps: the two first steps deal with $x\notin \mathcal{C}$,
the third with $x\in \mathcal{C}=\left( -\infty ,-M_{\tau }\right] \cup %
\left[ M_{\tau },+\infty \right) $.

\textbf{Part 1: }$x>0$

\textit{First step:} We split $\left[ 0,M_{\tau }\right] $ in two subsets, $%
\left[ 0,M_{\tau }\right] =\left[ 0,A_{\tau }\right] \cup \left[ A_{\tau
},M_{\tau }\right] $ with:  
\[
A_{\tau }=\omega /\left( M_{\tau }+1+\omega \right) 
\]
is chosen such that $0\leq x\leq A_{\tau }$ implies the following inequality
on the lower bound of the integral: $\left( \omega /x\right) -\lambda
-c>M_{\tau }$. Clearly $A_{\tau }\leq M_{\tau }$ because $\lambda >0>-\tau
-M_{\tau }$. Then: 
\[
\frac{1}{2c}\int_{\left( \omega /x\right) -\lambda -c}^{\left( \omega
/x\right) -\lambda +c}g\left( s\right) ds=\frac{C_{2}}{2c}\int_{\left(
\omega /x\right) -\lambda -c}^{\left( \omega /x\right) -\lambda +c}\log
^{2}\left\vert s\right\vert ds\leq C_{2}\log ^{2}\left\vert \left( \omega
/x\right) -\lambda +c\right\vert .
\]
Let: 
\[
\sup_{0\leq x\leq \omega /\left( c+M_{\tau }\right) }\sqrt{\left\vert
x\right\vert }\left( \log \left\vert \left( \omega /x\right) -\lambda
+c\right\vert \right) ^{2}=K_{1}\left( \omega ,c,\tau \right) <+\infty .
\]
The strictly positive $K_{1}\left( \omega ,c,\tau \right) $ exists because $%
\sqrt{x}\log ^{2}\left\vert \left( \omega /x\right) -\lambda +c\right\vert $
is bounded on $\left[ 0,A_{\tau }\right] $. The first condition reads: 
\[
\frac{1}{2c}\int_{\left( \omega /x\right) -\lambda -c}^{\left( \omega
/x\right) -\lambda +c}g\left( s\right) ds\leq \rho _{1}C_{1}/\sqrt{
\left\vert x\right\vert }, \qquad 0\leq x\leq A_{\tau }
\]
whenever 
\begin{equation}
C_{2}K_{1}\left( \omega ,c,\tau \right) \leq \rho _{1}C_{1}  \label{cond1}
\end{equation}
and $\rho _{1}$ will be fixed after the second step.\bigskip

\textit{Second step:} Now we turn to $A_{\tau }\leq x\leq M_{\tau }$. We
still have $g\left( x\right) =C_{1}/\sqrt{\left\vert x\right\vert }$ but we
need to focus on the bounds of the integral.

This time the lower bound of the integral $\left( \omega /x\right) -\lambda
-c\in \left[ -\lambda -\tau ,M_{\tau }\right] $ and the upper bound $\left(
\omega /x\right) -\lambda +c\in \left[ 2c-\lambda -\tau ,2c+M_{\tau }\right] 
$. We are going to require that $\left( \omega /x\right) -\lambda -c\geq
-M_{\tau }$ it suffices to take $\lambda +\tau \leq M_{\tau }$ and this
comes down to the following set of constraint on $\tau $: $\left\{ \tau \geq
c-\omega \right\} \cup \left\{ \tau \leq c-1\right\}.$ We keep the second
and assume once and for all that: 
\begin{equation}
\tau \leq c-1.  \label{constr-tau}
\end{equation}
Then for $x\in \left[ A_{\tau },M_{\tau }\right] $, 
\begin{align}
\frac{1}{2c}\int_{\left( \omega /x\right) -\lambda -c}^{\left( \omega
/x\right) -\lambda +c}g\left( s\right) ds& =\frac{1}{2c}\int_{\left( \omega
/x\right) -\lambda -c}^{M_{\tau }}g\left( s\right) ds+\frac{1}{2c}
\int_{M_{\tau }}^{\left( \omega /x\right) -\lambda +c}g\left( s\right) ds
\label{calc} \\
& \equiv \mathcal{I}_{1}+\mathcal{I}_{2}.  \nonumber
\end{align}

We want to make sure that the upper bound $\left( \omega /x\right) -\lambda
+c$ is larger than $M_{\tau }$. This will hold if $\left( \omega /M_{\tau
}\right) -\lambda +c\geq M_{\tau }$ hence if $2c-\tau -\lambda \geq M_{\tau
}.$ We imposed previously that $\lambda +\tau \leq M_{\tau }$. So: 
\[
\tau <c-\frac{\omega }{c}\Rightarrow M_{\tau }<c\Rightarrow 2c-\tau -\lambda
\geq M_{\tau }
\]
but the constraint $\tau <c-\omega /c$ is weaker than (\ref{constr-tau})
consequently (\ref{calc}) holds.

Focus on the first term $\mathcal{I}_{1}$ in (\ref{calc}) and consider: 
\[
\mathcal{I}_{1}=\frac{1}{2c}\int_{\left( \omega /x\right) -\lambda
-c}^{M_{\tau }}g\left( s\right) ds=\frac{1}{2c}\int_{\left( \omega /x\right)
-\lambda -c}^{M_{\tau }}\frac{C_{1}}{\sqrt{\left\vert s\right\vert }}ds.
\]
Consider the (only) two situations on the sign of $\left( \omega /x\right)
-\lambda -c=\left( \omega /x\right) -\left( 1+\omega \right) $.\newline
If $x<\omega /\left( 1+\omega \right)$ then $\left( \omega /x\right) -\left(
1+\omega \right) >0$ and $\mathcal{I}_{1}\leq \frac{C_{1}}{c}\sqrt{M_{\tau }}
$ Notice by the way and for further purpose that: 
\[
\sup_{x\in \left[ A_{\tau },\omega /\left( 1+\omega \right) \right] }\sqrt{
\left\vert x\right\vert }\mathcal{I}_{1}\leq \frac{C_{1}}{c}\sqrt{\frac{
\omega }{1+\omega }M_{\tau }}.
\]
If $x\geq \omega /\left( 1+\omega \right) $ then $\left( \omega /x\right)
-\left( 1+\omega \right) \leq 0$ and: 
\begin{align*}
\mathcal{I}_{1}& =\frac{1}{2c}\int_{\left( \omega /x\right) -\lambda
-c}^{M_{\tau }}g\left( s\right) ds=\frac{1}{2c}\int_{\left( \omega /x\right)
-\left( 1+\omega \right) }^{0}\frac{C_{1}}{\sqrt{\left\vert s\right\vert }}
ds+\frac{1}{2c}\int_{0}^{M_{\tau }}\frac{C_{1}}{\sqrt{\left\vert
s\right\vert }}ds \\
& =\frac{C_{1}}{c}\left[ \sqrt{\left\vert \left( \omega /x\right) -\left(
1+\omega \right) \right\vert }+\sqrt{M_{\tau }}\right].
\end{align*}
Again: 
\[
\sqrt{\left\vert x\right\vert }\mathcal{I}_{1}\leq \frac{C_{1}}{c}\left[ 
\sqrt{\left\vert x\left( 1+\omega \right) -\omega \right\vert }+\sqrt{
xM_{\tau }}\right].
\]

From the bounds above we see that: 
\[
\sup_{x\in \left[ A_{\tau },M_{\tau }\right] }\sqrt{\left\vert x\right\vert }
\mathcal{I}_{1}\leq \frac{C_{1}}{c}\left[ \sqrt{\left\vert M_{\tau }\left(
1+\omega \right) -\omega \right\vert }+M_{\tau }\right].
\]
The reader will soon understand why we need to make sure that the right had
side in equation above is strictly under $C_{1}$. It is not hard to see that
the function $\tau \longmapsto \sqrt{\left\vert M_{\tau }\left( 1+\omega
\right) -\omega \right\vert }+M_{\tau }$ is increasing and continuous on $%
\left[ 0,c-1\right] $. If we prove that for some $\delta \in \left] 0,1 %
\right[ $: 
\[
\frac{1}{c}\left[ \sqrt{\left\vert M_{0}\left( 1+\omega \right) -\omega
\right\vert }+M_{0}\right] =1-3\delta <1,
\]

then the existence of some $\tau ^{+}>0$ such that: 
\begin{equation}
\frac{1}{c}\left[ \sqrt{\left\vert M_{\tau ^{+}}\left( 1+\omega \right)
-\omega \right\vert }+M_{\tau ^{+}}\right] =1-2\delta <1  \label{delta}
\end{equation}
will be granted.\ But $\frac{1}{c}\left[ \sqrt{\left\vert M_{0}\left(
1+\omega \right) -\omega \right\vert }+M_{0}\right] =\frac{1}{c}\left[ \sqrt{
\frac{\omega }{c}\lambda }+\frac{\omega }{c}\right].$ \newline
If we assume that $\lambda <\omega /c<\left( 1+c\right) /4$ (assumption $%
\mathbf{B}_{1}$) then since $c>1$: 
\[
\frac{1}{c}\left[ \sqrt{\frac{\omega }{c}\lambda }+\frac{\omega }{c}\right]
< \frac{1}{2}\left( 1+\frac{1}{c}\right) <1.
\]

We turn to $\mathcal{I}_{2}$ in (\ref{calc}): 
\[
\mathcal{I}_{2}=\frac{1}{2c}\int_{M_{\tau }}^{\left( \omega /x\right)
-\lambda +c}g\left( s\right) ds=\frac{C_{2}}{2c}\int_{M_{\tau }}^{\left(
\omega /x\right) -\lambda +c}\left( \log s\right) ^{2} ds\leq \frac{C_{2}}{ 
\sqrt{\left\vert x\right\vert }}K_{2}\left( \omega ,c,\tau ^{+}\right),
\]
where: 
\[
K_{2}\left( \omega ,c,\tau ^{+}\right) =\sup_{x\in \left[ A_{\tau },M_{\tau
} \right] }\frac{\sqrt{\left\vert x\right\vert }}{2c}\int_{M_{\tau
}}^{\left( \omega /x\right) -\lambda +c}\left( \log s\right) ^{2}ds.
\]

Set finally $\rho _{1}^{+}=1-\delta <1$.

From (\ref{calc}) we get: 
\begin{align*}
\frac{1}{2c}\int_{\left( \omega /x\right) -\lambda -c}^{\left( \omega
/x\right) -\lambda +c}g\left( s\right) ds& \leq \frac{C_{1}}{\sqrt{
\left\vert x\right\vert }}.\left( 1-2\delta \right) +\frac{C_{2}}{\sqrt{
\left\vert x\right\vert }}K_{2}\left( \omega ,c,\tau ^{+}\right) \\
& \leq \rho _{1}^{+}\frac{C_{1}}{\sqrt{\left\vert x\right\vert }},
\end{align*}
whenever holds the new condition: 
\begin{equation}
C_{2}K_{2}\left( \omega ,c,\tau \right) \leq C_{1}\delta.  \label{cond2}
\end{equation}
Finally comparing (\ref{cond1}) and (\ref{cond2}), we see that both
conditions cannot be incompatible. Accurate choices of the couple $\left(
C_{1}^{+},C_{2}^{+}\right) $ are given by the summary bound: 
\begin{equation}
C_{2}^{+}\leq C_{1}^{+}\min \left( \frac{\delta }{K_{2}},\frac{1-\delta }{
K_{1}}\right).  \label{condPart1}
\end{equation}
It is now basic to see that the quadruple $\left( C_{1}^{+},C_{2}^{+},\tau
^{+},\rho _{1}^{+}\right) $ yields the drift condition (\ref{drift_condition}%
) for $x\notin \mathcal{C}$.\bigskip

\textit{Third step:} The remaining step is to check the inequality for some $%
\rho _{2}$: 
\[
\int_{\mathbb{R}}g\left( t\right) P\left( t,x\right) dt\leq \rho
_{1}^{+}g\left( x\right) +\rho _{2},
\]
for any $x$ in $\mathcal{C}$ -that is any $\left\vert x\right\vert >M_{\tau
} $ (rather $x>M_{\tau }$ here as explained above since $x>0$). We see that: 
\[
0\leq \frac{\omega }{x}\leq \frac{\omega }{M_{\tau }},
\]
and: 
\[
\frac{1}{2c}\int_{\left( \omega /x\right) -\lambda -c}^{\left( \omega
/x\right) -\lambda +c}g\left( s\right) ds\leq \frac{1}{2c}\int_{-\lambda
-c}^{2c-\tau -\lambda }g\left( s\right) ds\leq \frac{1}{2c}\int_{-\left(
1+\omega \right) }^{3c-\left( 1+\omega \right) }g\left( s\right) ds.
\]
The values of the constants $C_{1}$ and $C_{2}$ were fixed above. Then
denote: 
\[
\rho _{2}^{+}=\frac{1}{2c}\int_{-\left( 1+\omega \right) }^{3c-\left(
1+\omega \right) }g\left( s\right) ds>0,
\]
then clearly for any $x$ in $\mathcal{C}$: 
\[
\int_{\mathbb{R}}g\left( t\right) P\left( t,x\right) dt\leq \rho _{2}^{+},
\]
so that (\ref{drift_condition}) holds.\bigskip

\textbf{Part 2 (}$x\leq 0$\textbf{)}

We go on with $x<0$ and $\lambda >0,$ set $y=-x\geq 0,$ 
\[
\int_{\mathbb{R}}g\left( t\right) P\left( t,x\right) dt=\frac{1}{2c}\int_{ 
\mathbb{\lambda -}\left( \omega /x\right) -c}^{\mathbb{\lambda -}\left(
\omega /x\right) +c}g\left( s\right) ds=\frac{1}{2c}\int_{\left( \omega
/y\right) +\lambda -c}^{\left( \omega /y\right) +\lambda +c}g\left( s\right)
ds.
\]
Since $g$ is even and in view of the proposed $\mathcal{C}$ we just have to
prove exactly the following drift condition with $x>0$: 
\[
\frac{1}{2c}\int_{\left( \omega /x\right) +\lambda -c}^{\left( \omega
/x\right) +\lambda +c}g\left( s\right) ds\leq \rho _{1}g\left( x\right)
+\rho _{2}1_{x\in \mathcal{C}}.
\]
\bigskip \textit{First step:} Take $x\notin \mathcal{C}$. We split $\left[
0,M_{\tau }\right] $ in two subsets, $\left[ 0,M_{\tau }\right] =\left[
0,B_{\tau }\right] \cup \left[ B_{\tau },M_{\tau }\right] $ with $B_{\tau
}=\omega /\left( M_{\tau }-\lambda +c\right) $ is chosen such that $0\leq
x\leq B_{\tau }$ implies the following inequality on the lower bound of the
integral: $\left( \omega /x\right) +\lambda -c>M_{\tau }$. Clearly $B_{\tau
}\leq M_{\tau }$ for all $\tau $. Then: 
\[
\frac{1}{2c}\int_{\left( \omega /x\right) +\lambda -c}^{\left( \omega
/x\right) +\lambda +c}g\left( s\right) ds=\frac{C_{2}}{2c}\int_{\left(
\omega /x\right) +\lambda -c}^{\left( \omega /x\right) +\lambda +c}\log
^{2}\left\vert s\right\vert ds\leq C_{2}\log ^{2}\left\vert \left( \omega
/x\right) +\lambda +c\right\vert .
\]
Let: 
\[
\sup_{0\leq x\leq B_{\tau }}\sqrt{\left\vert x\right\vert }\log
^{2}\left\vert \left( \omega /x\right) +\lambda +c\right\vert =K_{1}\left(
\omega ,c,\tau \right) <+\infty .
\]
The strictly positive $K_{1}\left( \omega ,c,\tau \right) $ exists because $%
\sqrt{\left\vert x\right\vert }\log ^{2}\left\vert \left( \omega /x\right)
+\lambda +c\right\vert $ is bounded on $\left[ 0,B_{\tau }\right] $. The
initial condition reads: 
\[
\frac{1}{2c}\int_{\left( \omega /x\right) +\lambda -c}^{\left( \omega
/x\right) +\lambda +c}g\left( s\right) ds\leq \rho _{1}C_{1}/\sqrt{
\left\vert x\right\vert },\qquad 0\leq x\leq B_{\tau }
\]
whenever $C_{2}K_{1}\left( \omega ,c,\tau \right) \leq \rho _{1}C_{1}$ and $%
\rho _{1}$ will be fixed later.\bigskip

\textit{Second step:} Now we turn to $B_{\tau }\leq x\leq M_{\tau }$. We
still have $g\left( x\right) =C_{1}/\sqrt{\left\vert x\right\vert }$ but we
need to focus on the bounds of the integral.

This time the lower bound of the integral $\left( \omega /x\right) +\lambda
-c\in \left[ \lambda -\tau ,M_{\tau }\right] $ and the upper bound $\left(
\omega /x\right) +\lambda +c\in \left[ \lambda +2c-\tau ,M_{\tau }+2c\right] 
$. If we assume that $\tau \leq \lambda $ then $\left( \omega /x\right)
+\lambda -c\geq 0\geq -M_{\tau }.$ Besides in order that $\lambda +2c-\tau
\geq M_{\tau }$ we just have need that $2c\geq M_{\tau }$ or $\tau \leq
c-\omega /\left( 2c\right).$ As a consequence the assumption: 
\begin{equation}
\tau \leq \min \left( \lambda ,c-\frac{\omega }{2c}\right)  \label{condthau2}
\end{equation}
allows to write for $x\in \left[ B_{\tau },M_{\tau }\right] $, 
\begin{align*}
\frac{1}{2c}\int_{\left( \omega /x\right) +\lambda -c}^{\left( \omega
/x\right) +\lambda +c}g\left( s\right) ds& =\frac{1}{2c}\int_{\left( \omega
/x\right) +\lambda -c}^{M_{\tau }}g\left( s\right) ds+\frac{1}{2c}
\int_{M_{\tau }}^{\left( \omega /x\right) +\lambda +c}g\left( s\right) ds \\
& \equiv \mathcal{I}_{1}+\mathcal{I}_{2},
\end{align*}

with non-null $\mathcal{I}_{1}$ and $\mathcal{I}_{2}$. Focus on: 
\begin{eqnarray*}
\mathcal{I}_{1} &=&\frac{C_{1}}{2c}\int_{\left( \omega /x\right) +\lambda
-c}^{M_{\tau }}\frac{1}{\sqrt{s}}ds=\frac{C_{1}}{c}\left[ \sqrt{M_{\tau }}- 
\sqrt{\left( \omega /x\right) +\lambda -c}\right] \\
\sup_{x\in \left[ B_{\tau },M_{\tau }\right] }\sqrt{x}\mathcal{I}_{1} &\leq
& \frac{C_{1}}{c}M_{\tau }.
\end{eqnarray*}
At last we see that for $M_{\tau }\leq 1$ i.e. $\tau \leq c-\omega $, $%
\sup_{x\in \left[ B_{\tau },M_{\tau }\right] }\sqrt{x}\mathcal{I}_{1}\leq
C_{1}/c$. This condition combined with (\ref{condthau2}) let us set in the
sequel: 
\[
\tau \leq \tau ^{-}=\min \left( \lambda ,c-\omega \right).
\]
We turn to $\mathcal{I}_{2}$: 
\[
\frac{1}{2c}\int_{M_{\tau }}^{\left( \omega /x\right) +\lambda +c}g\left(
s\right) ds=\frac{C_{2}}{2c}\int_{M_{\tau }}^{\left( \omega /x\right)
+\lambda +c}\log ^{2}\left\vert s\right\vert ds\leq \frac{C_{2}}{\sqrt{x}}
K_{2}\left( \omega ,c,\tau ^{-}\right),
\]
with: 
\[
K_{2}^{-}=K_{2}\left( \omega ,c,\tau ^{-}\right) =\sup_{x\in \left[ B_{\tau
},M_{\tau }\right] }\frac{\sqrt{\left\vert x\right\vert }}{2c}\int_{M_{\tau
}}^{\left( \omega /x\right) +\lambda +c}\left( \log s\right) ^{2}ds.
\]

Set finally $\rho _{1}^{-}=\left( 1+1/c\right) /2<1$. From all that was done
above we get: 
\[
\int_{\mathbb{R}}g\left( t\right) P\left( t,x\right) dt\leq \frac{C_{1}}{ 
\sqrt{\left\vert x\right\vert }}\frac{1}{c}+\frac{C_{2}}{\sqrt{\left\vert
x\right\vert }}K_{2}^{-}\leq \rho _{1}^{-}\frac{C_{1}}{\sqrt{\left\vert
x\right\vert }},
\]
whenever: 
\[
C_{2}K_{2}\left( \omega ,c,\tau ^{-}\right) \leq \frac{1-1/c}{2}C_{1}.
\]
This will be combined with the constraint of the first step $%
C_{2}K_{1}^{-}\leq \rho _{1}C_{1}$ (we denoted $K_{1}\left( \omega ,c,\tau
^{-}\right) =K_{1}^{-}$). The new condition: 
\begin{equation}
C_{2}^{-}\leq C_{1}^{-}\min \left( \frac{\rho _{1}^{-}}{K_{1}^{-}},\frac{
1-1/c}{2K_{2}^{-}}\right)  \label{condPart2}
\end{equation}
ensures that 
\[
\frac{1}{2c}\int_{\left( \omega /x\right) +\lambda -c}^{\left( \omega
/x\right) +\lambda +c}g\left( s\right) ds\leq \rho _{1}^{-}g\left( x\right) 
\text{ for } x\in \left[ 0,M_{\tau }\right]. 
\]%
\bigskip

\textit{Third step:} The remaining step is to check the inequality: 
\[
\int_{\mathbb{R}}g\left( t\right) P\left( t,x\right) dt\leq \rho
_{1}^{-}g\left( x\right) +\rho _{2},
\]
for any $x$ in $\mathcal{C}$ -that is here any $x>M_{\tau }$. Adapting the
method given above is straightforward and leads to the desired result with a
given $\rho _{2}^{-}$.\bigskip

We are ready to conclude. Take 
\[
C_{2}^{\ast }= C_{1}^{\ast }\min \left( \frac{\delta }{K_{2}^{+}},\frac{
1-\delta }{K_{1}^{+}},\frac{\rho _{1}^{-}}{K_{1}^{-}},\frac{1-1/c}{%
2K_{2}^{-} }\right) .
\]
Conditions (\ref{condPart1}) and (\ref{condPart2}) hold for the couple $%
\left( C_{1}^{\ast },C_{2}^{\ast }\right).$ For such a couple we have: 
\begin{eqnarray*}
\int_{\mathbb{R}}g\left( t\right) P\left( t,x\right) dt &\leq &\rho
_{1}^{+}g\left( x\right) +\rho _{2}^{+},\quad x>0, \\
\int_{\mathbb{R}}g\left( t\right) P\left( t,x\right) dt &\leq &\rho
_{1}^{-}g\left( x\right) +\rho _{2}^{-},\quad x\leq 0,
\end{eqnarray*}
and for all $x$: 
\[
\int_{\mathbb{R}}g\left( t\right) P\left( t,x\right) dt\leq \max \left( \rho
_{1}^{+},\rho _{1}^{-}\right) g\left( x\right) +\max \left( \rho
_{2}^{+},\rho _{2}^{-}\right).
\]

This finishes the proof of the Lemma.

\end{document}